\documentclass[12pt]{article}
\usepackage{amsmath,amssymb}

\newcommand{\ncm}{\newcommand}
\ncm{\Inn}{\mbox{\rm Inn}} \ncm{\Ap}{\mbox{$\overline{\rm Inn}$}} 
\ncm{\Ext}{\mbox{\rm Ext}} \ncm{\Ex}{\mbox{\rm Ex}} 
\ncm{\OExt}{\mbox{\rm OrderExt}} \ncm{\AI}{\mbox{\rm AInn}} 
\ncm{\HI}{\mbox{\rm HInn($A$)}} \ncm{\Aut}{\mbox{\rm Aut}} 
\ncm{\Mal}{\mbox{$M_{\alpha}$}} \ncm{\Aff}{\mbox{${\rm Aff}$}} 
\ncm{\id}{\mbox{\rm id}} \ncm{\Ker}{\mbox{\rm Ker}} 
\ncm{\BE}{\begin{eqnarray*}} \ncm{\EE}{\end{eqnarray*}} 
\ncm{\lra}{\mbox{$\longrightarrow$}} \ncm{\Hom}{\mbox{\rm Hom}} 
\ncm{\calU}{{\cal U}} \ncm{\el}{\ell} 
 \ncm{\ad}{{\rm ad}} 
 \ncm{\diag}{{\rm diag}}    
 \ncm{\rank}{{\rm rank}}       
 \ncm{\Tr}{{\rm Tr}} 
 \ncm{\Arg}{{\rm Arg}} 
\ncm{\Alg}{\mbox{\rm Alg}} \ncm{\Conv}{\mbox{\rm Conv}} 
 \ncm{\D}{{ D}} \ncm{\calE}{{\mathcal{E}} }    \ncm{\calJ}{{\mathcal{J}} }     
\ncm{\calS}{{\mathcal{S}} } \ncm{\UHF}{{\rm UHF}}      

\ncm{\cstar}{C$^{*}$-algebra} \ncm{\cstars}{C$^{*}$-algebras} 
\ncm{\subcstar}{C$^{*}$-subalgebra}  
\ncm{\ra}{\mbox{$\rightarrow$}} \ncm{\la}{\mbox{$\leftarrow$}} 
\ncm{\hra}{\hookrightarrow} \ncm{\da}{\mbox{$\downarrow$}} 
\ncm{\se}{\mbox{$\searrow$}} \ncm{\al}{\mbox{$\alpha $}} 
\ncm{\del}{\mbox{$\delta$}} \ncm{\supp}{\mbox{\rm supp}} 
\ncm{\Ad}{\mbox{\rm Ad}} \ncm{\CAR}{\mbox{$M_{2^{\infty}}$}} 
\ncm{\ep}{\mbox{$\epsilon > 0$}} \ncm{\ol}{\overline} 
\ncm{\Mninf}{\mbox{$M_{n^{\infty}}$}} \ncm{\MR}{M. R\o{}rdam} 
\ncm{\Range}{\mbox{\rm Range}} 
\ncm{\vo}{}%{\bf}
\ncm{\ch}{}%{\it}
\ncm{\CMP}{Comm. Math. Phys.} \ncm{\add}{} 
\ncm{\tilsig}{\tilde{\sigma}} \ncm{\dist}{{\rm 
dist}}\ncm{\eps}{\epsilon}  \ncm{\calL}{{\mathcal{L}}}   
 \ncm{\E}{{\mathcal{E}} }    \ncm{\M}{{\mathcal{M}} }
 \ncm{\calO}{{\mathcal{O}}} \ncm{\F}{{\mathcal{F}}} \ncm{\G}{{\mathcal{G}} }  
\ncm{\calH}{{\mathcal{H}}}   \ncm{\calC}{{\mathcal{C}}} 
\ncm{\lan}{{\langle}}\ncm{\ran}{{\rangle}}   
\ncm{\calK}{{\mathcal{K}}}    \ncm{\Spec}{{\rm Spec}}  
\ncm{\calP}{{\mathcal{P}}} \ncm{\Hil}{{\mathcal{H}}} 
\ncm{\U}{{\mathcal{U}}}                \ncm{\A}{\mathcal{A}}    
\ncm{\Nr}{\mathcal{N}}  %\ncm{\calC}{\mathcal{C}}     

\newtheorem{theo}{Theorem}[section]
\newtheorem{cor}[theo]{Corollary}
\newtheorem{lem}[theo]{Lemma}
\newtheorem{prop}[theo]{Proposition}
\newtheorem{remark}[theo]{Remark}
\newtheorem{definition}[theo]{Definition}
\newtheorem{example}[theo]{Example}
\newtheorem{property}[theo]{Property}

\newenvironment{rem}{\begin{remark} \rm}{\end{remark}}
\newenvironment{pf}{{\it Proof.}}{\hfill$\square$\vspace{3mm}}

\ncm{\R}{\mbox{\bf R}} \ncm{\Z}{\mbox{\bf Z}} \ncm{\T}{{\bf T}} 
\ncm{\TT}{\T$^{2}$} \ncm{\N}{\mbox{\bf N}} \ncm{\C}{\mbox{\bf C}}

\title{Rohlin flows on the Cuntz algebra $\calO_\infty$}
\bigskip
\author{Ola Bratteli \\
{\small Department of Mathematics, University of Oslo}\\ {\small 
Blindern, P.O.Box 1053, N-0316, Norway} 
\smallskip\\ Akitaka Kishimoto \\ {\small 
Department of Mathematics, Hokkaido University, Sapporo 060-0810, 
Japan}\\ and\\ 
 Derek W. Robinson\\
{\small Centre for Mathematics and its Applications, Australian 
National University}\\ {\small Canberra, ACT 0200, Australia}} 
\date{\small January, 2005}
\begin{document}
\maketitle 

\begin{abstract} 
It is shown that certain quasi-free flows on the Cuntz algebra 
$\calO_\infty$ have the Rohlin property and therefore are 
cocycle-conjugate with each other. This, in particular, shows that 
any unital separable nuclear purely infinite simple \cstar\ has a 
Rohlin flow. 
\end{abstract}                  
                    
\section{Introduction}    
We are concerned here with Rohlin flows; a flow $\alpha$ on a 
unital \cstar\ $A$ is said to have the Rohlin property (or to be a 
Rohlin flow) if for any $p\in\R$ there is a central sequence 
$(u_n)$ in $\U(A)$, the unitary group of $A$, such that 
$\max_{|t|\leq1}\|\alpha_t(u_n)-e^{ipt}u_n\|\ra0$ as $n\ra\infty$. 
A major consequence of this property can be paraphrased as {\em 
any $\alpha$-cocycle is almost a coboundary}. This consequence, 
combined with enough information on $\U(A)$, may lead us to a 
classification theory of Rohlin flows up to cocycle conjugacy. 
This is a goal we have in mind (see \cite{K96b,K02,K03,Kp,K05}). 

Since the property is rather stringent, it is not easy to present 
a Rohlin flow in general. But we managed to give Rohlin flows on 
the Cuntz algebra $\calO_n$ with $n$ finite; moreover we can 
identify the quasi-free flows which have the Rohlin property. In 
this paper we show that certain quasi-free flows on $\calO_\infty$ 
have the Rohlin property. Hence it follows that any unital \cstar\ 
$A$ with $A\cong A\otimes \calO_\infty$ has Rohlin flows; the 
class of such $A$ includes all unital separable nuclear purely 
infinite simple \cstars, due to Kirchberg.

We have left some quasi-free flows $\calO_\infty$ undecided 
whether they have the Rohlin property or not. But, as in 
\cite{K02}, we show that all the Rohlin flows on $\calO_\infty$ 
are cocycle conjugate with each other in the class of quasi-free 
flows. This is true for a wider class of flows. Noting that there 
is a certain maximal abelian C$^*$-subalgebra $\calC_\infty$ of 
$\calO_\infty$ whose elements the quasi-free flows fix, we show 
that Rohlin flows are cocycle-conjugate in the class of flows 
which are $C^{1+\eps}$ on $\calC_\infty$ (see below for details 
and note that our terminology of quasi-free flows is restrictive). 

We will now describe the contents more precisely. 

For each $n=2,3,\ldots$ the Cuntz algebra $\calO_n$ is generated 
by $n$ isometries $s_1,s_2,\ldots,s_n$ such that 
$\sum_{k=1}^ns_ks_k^*=1$. For $n=\infty$ the Cuntz algebra 
$\calO_\infty$ is generated by a sequence $(s_1,s_2,\ldots)$ of 
isometries such that $\sum_{k=1}^ns_ks_k^*\leq 1$ for all $n$. It 
is shown in \cite{Cun1} that $\calO_n$ with $n=2,3,\ldots$ or 
$n=\infty$  is a simple purely infinite nuclear \cstar. 

For a finite (resp. infinite) sequence $(p_1,p_2,\ldots,p_n)$ in 
$\R$ we define a flow $\alpha$, called a quasi-free flow, on 
$\calO_n$ (resp. $\calO_\infty$) by 
 $$
 \alpha_t(s_k)=e^{ip_kt}s_k.
 $$                                                     
(In the case $n=\infty$, a more general flow can be induced by a 
unitary flow $U$ on the closed linear subspace $\Hil$ spanned by 
$s_1,s_2,\ldots$, where the inner product $\lan\cdot\,,\cdot\ran$ 
is given by $y^*x=\lan x,y\ran 1,\ x,y\in \Hil$, if the generator 
of $U$ is not diagonal. But we will exclude them from the 
quasi-free flows in this paper.) It is known in \cite{K80,KK} that 
if $p_1,p_2,\ldots$ generate $\R$ as a closed subsemigroup, then 
the crossed product $\calO_n\times_\alpha\R$ is simple and purely 
infinite (whether $n$ is finite or infinite). It is also known in 
\cite{K02,K05} that if $n$ is finite, the flow $\alpha$ has the 
Rohlin property if and only if $\calO_\infty\times_\alpha\R$ is 
simple and purely infinite. For $n=\infty$, it is known in 
\cite{K80,K96b} that if $\alpha$ has the Rohlin property then 
$\calO_\infty\times_\alpha\R$ is simple and purely infinite. In 
this paper we shall give a partial converse to this fact: 

\begin{theo}
Let $(p_k)$ be an infinite sequence in $\R$ such that 
$p_1,p_2,\ldots,p_n$ generate $\R$ as a closed subsemigroup for 
some $n$. Then the quasi-free flow $\alpha$ on $\calO_\infty$ 
defined by $\alpha_t(s_k)=e^{ip_kt}s_k$ has the Rohlin property. 
\end{theo}     

We shall prove that each $\alpha_t$ is $\alpha$-invariantly 
approximately inner, i.e., for each $t\in\R$ there is a sequence 
$(u_n)$ in $\U(\calO_\infty)$ such that 
$\alpha_t(x)=\lim\Ad\,u_n(x),\ x\in\calO_\infty$ and 
$\max_{s\in[0,1]}\|\alpha_s(u_n)-u_n\|\ra0$. Then we would get the 
above theorem, by \cite{K03,Kp}, from the fact that 
$\calO_\infty\times_\alpha\R$ is simple and purely infinite.  

Let $\calE_n$ be the C$^*$-subalgebra of 
$\calO_\infty=C^*(s_1,s_2,\ldots)$ generated by 
$s_1,s_2,\ldots,s_n$. Then $\calE_n$ is left invariant under 
$\alpha$ and the union $\bigcup_n\calE_n$ is dense in 
$\calO_\infty$. Hence, to prove the assertion in the previous 
paragraph, it suffices to show that $\alpha|\calE_n$ is 
$\alpha$-invariantly approximately inner for all large $n$. Let us 
state formally:

\begin{prop}\label{En}
Let $s_1,s_2,\ldots,s_n$ be isometries such that
 $$
 \sum_{k=1}^ns_ks_k^*\lneqq 1
 $$
and let $\calE_n$ be the \cstar\ generated by these 
$s_1,\ldots,s_n$. Let $(p_1,p_2,\ldots,p_n)$ be a finite sequence 
in $\R$ such that $p_1,\ldots,p_n$ generate $\R$ as a closed 
subsemigroup and define a quasi-free flow $\alpha$ on $\calE_n$ by 
$\alpha_t(s_k)=e^{ip_kt}s_k$. Then each $\alpha_t$ is 
$\alpha$-invariantly approximately inner. \end{prop}

To prove this we use the following facts. Let $\calJ_n$ be the 
ideal of $\calE_n$ generated by $e_n^0=1-\sum_{k=1}^ns_ks_k^*$. 
Then $\calJ_n$ is isomorphic to the C$^*$-algebra $\calK$ of 
compact operators (on a separable infinite-dimensional Hilbert 
space) and is left invariant under $\alpha$. The quotient 
$\calE_n/\calJ_n$ is isomorphic to $\calO_n$ by mapping 
$s_k+\calJ_n$ into $s_k$ (the latter $s_k$'s satisfy the equality 
$\sum_{k=1}^ns_ks_k^*=1$ and generate $\calO_n$). By the 
assumption on $(p_k)$ the induced flow $\dot{\alpha}$ on $\calO_n$ 
has the Rohlin property \cite{K05}, from which follows that each 
$\dot{\alpha}_t$ is $\alpha$-invariantly approximately inner. We 
will translate this property to $\alpha_t$ on $\calE_n$ by using 
the fact that $\calJ_n\cong \calK$. See Section 3 for details.

Before embarking on the proof of the above proposition, we will 
have to prove that if $\alpha$ is a Rohlin flow on $\calO_n$, then 
each $\alpha_t$ is not only $\alpha$-invariantly approximately 
inner but also $\alpha$-invariantly asymptotically inner, i.e., 
there is a continuous map $u:[0,\infty)\ra\U(\calO_n)$ such that 
$\alpha_t(x)=\lim_{s\rightarrow \infty}\Ad\,u(s)(x)$ for 
$x\in\calO_n$ and $\max_{t_1\in 
[0,1]}\|\alpha_{t_1}(u(s))-u(s)\|\ra0$. This will be proved for a 
wider class of \cstars\ (see \ref{A1} for details). (As a matter 
of fact we do not know of a single example of $\alpha$ without the 
above property of $\alpha$-invariant asymptotical innerness if it 
has covariant irreducible representations; we expect that this 
property holds fairly in general whether it has the Rohlin 
property or not.)  

As a corollary to the above theorem we get that any purely 
infinite simple separable nuclear \cstar\ has a Rohlin flow; 
because such a \cstar\ $A$ satisfies that $A\cong A\otimes 
\calO_\infty$ due to Kirchberg (see \cite{KP1}) and a flow 
$\alpha$ on $\calO_\infty$ induces a flow on $A$ via 
$\id\otimes\alpha$ on $A\otimes\calO_\infty$ which has the Rohlin 
property if $\alpha$ has. 

Let $\calC_\infty$ denote the C$^*$-subalgebra of $\calO_\infty$ 
generated by $s_{i_1}s_{i_2}\cdots s_{i_k}s_{i_k}^*\cdots 
s_{i_1}^*$ with all finite sequences $(i_1,i_2,\ldots,i_k)$ in 
$\N$. Then $\calC_\infty$ is a weakly regular maximal abelian 
C$^*$-subalgebra of $\calO_\infty$ ({\em weakly regular} in the 
sense that $\{u\in\mathcal{PI}(\calO_\infty)\ |\ uu^*,u^*u\in 
\calC_\infty,\  u\calC_\infty u^*=\calC_\infty uu^*\}$ generates 
$\calO_\infty$, where $\mathcal{PI}(\calO_\infty)$ is the set of 
partial isometries of $\calO_\infty$). Moreover there is a 
projection of norm one of $\calO_\infty$ onto $\calC_\infty$ and 
there is a character of $\calC_\infty$ which extends uniquely to a 
state of $\calO_\infty$. (When a weakly regular masa satisfies 
these two additional conditions, we will say that it is a {\em 
weak Cartan} masa.) We note that if $\alpha$ is a quasi-free flow 
(in our sense) then $\alpha_t$ is the identity on $\calC_\infty$; 
in other words, if $\delta_\alpha$ denotes the generator of 
$\alpha$, then $D(\delta_\alpha)\supset \calC_\infty$ and 
$\delta_\alpha|\calC_\infty=0$. We consider the following 
condition for a flow $\gamma$ on $\calO_\infty$: 
$D(\delta_\gamma)\supset \calC_\infty$ and $\sup_{x\in 
\calC_\infty,\, \|x\|\leq1} \|(\gamma_t-\id)\delta_\gamma(x)\|$ 
converges to zero as $t\ra0$; which we express by saying that {\em 
$\gamma$ is $C^{1+\eps}$ on $\calC_\infty$} below. This is 
obviously satisfied if $\gamma$ is $C^2$ on $\calC_\infty$ or 
$D(\delta_\gamma^2)\supset \calC_\infty$ (because then 
$\delta_\alpha^2|\calC_\infty$ is bounded).   

We will also show:

\begin{cor} \label{CC}
Any two Rohlin flows on $\calO_\infty$ are cocycle conjugate with 
each other if they are $C^{1+\eps}$ on $\calC_\infty$. 
\end{cor}            
                         
The proof consists of two parts. In the first part we show that if 
the flow $\gamma$ is $C^{1+\eps}$ on $\calC_\infty$ then 
$\delta_\gamma|\calC_\infty$ is inner, i.e., there is an 
$h=h^*\in\calO_\infty$ such that $\delta_\gamma(x)=\ad\,ih(x),\ 
x\in\calC_\infty$ (see \ref{oinf}). Thus we can assume, by inner 
perturbation, that $\delta_\gamma|\calC_\infty=0$. In the second 
part we show that any two Rohlin flows are cocycle-conjugate with 
each other if they fix each element of $\calC_\infty$ (see 
\ref{M}).  

{\small {\em Acknowledgement.} One of the authors (A.K.) visited 
at Australian National University in March, 2004 and at University 
of Oslo in August-September, 2004 during this collaboration. He 
acknowledges partial financial supports from these institutions.}     
                                                           
\section{Rohlin property} 

In this section we consider the class of purely infinite simple 
nuclear separable \cstar\ satisfying the universal coefficient 
theorem, which is classified by Kirchberg and Phillips 
\cite{KP1,KP2} in terms of K-theory.  

Let $A$ be a unital \cstar\ of the above class. Let 
$\ell^\infty(A)$ be the \cstar\ of bounded sequences in $A$ and 
let, for a free ultrafilter $\omega$ on $\N$, $c_\omega(A)$ be the 
ideal of $\ell^\infty(A)$ consisting  of $x=(x_n)$ with 
$\lim_{\omega}\|x_n\|=0$. If $\alpha$ is a flow on $A$, i.e., a 
strongly continuous one-parameter automorphism group of $A$, we 
can define an action of $\R$ on $\ell^\infty(A)$ by $t\mapsto 
(\alpha_t(x_n))$ for $x=(x_n)$. Let $\ell_\alpha^\infty(A)$ be the 
maximal C$^*$-subalgebra of $\ell^\infty(A)$ on which this action 
is continuous; we will denote this flow by $\alpha$. We set 
 $$
 A^\omega=\ell^\infty(A)/c_\omega(A),\ \ \ A_\alpha^\omega=\ell_\alpha^\omega(A)/c_\omega(A).
 $$
We embed $A$ into $\ell^\infty_\alpha(A)$ by constant sequences. 
Since $A\cap c_\omega(A)=\{0\}$, we regard $A$ as a 
C$^*$-subalgebra of $A_\alpha^\omega\subset A^\omega$.   

We recall the following result \cite{K03,Kp}: 

\begin{theo}\label{X}  
Let $A$ be a unital separable nuclear purely infinite simple 
\cstar\ satisfying the universal coefficient theorem and let 
$\alpha$ be a flow on $A$. Then the following conditions are 
equivalent. 
 \begin{enumerate}
 \item $\alpha$ has the Rohlin property.
 \item $(A'\cap A^\omega_\alpha)^\alpha$ is purely infinite and 
 simple, $K_0((A'\cap 
 A_\alpha^\omega)^\alpha)\cong K_0(A'\cap A^\omega)$ induced by the embedding,
  and $\Spec(\alpha|A'\cap A^\omega_\alpha)=\R$.
 \item The crossed product $A\times_\alpha\R$ is purely infinite 
 and simple and the dual action $\hat{\alpha}$ has the Rohlin 
 property.
 \item The crossed product $A\times_{\alpha}\R$ is purely infinite and 
 simple and $\alpha_{t_0}$ 
 is $\alpha$-invariantly approximately inner for every $t_0\in\R$.
 \end{enumerate}    
If the above conditions are satisfied, it also follows that  
$K_1((A'\cap A_\alpha^\omega)^\alpha)\cong K_1(A'\cap A^\omega)$, 
which is induced by the embedding. 
\end{theo} 

In the last condition of the above theorem, $\alpha_{t_0}$ (for a 
fixed $t_0$) is $\alpha$-invariantly approximately inner if there 
is a sequence $(u_n)$ in $\U(A)$ such that 
$\alpha_{t_0}=\lim\Ad\,u_n$ and 
$\max_{t\in[0,1]}\|\alpha_t(u_n)-u_n\|\ra0$. We will strengthen 
this condition as follows.  

\begin{lem} \label{A1} 
Let $\alpha$ be a Rohlin flow on a unital \cstar\ $A$ of the above 
class (or in particular $\calO_n$). Then each $\alpha_{t_0}$ is 
$\alpha$-invariantly asymptotically inner, i.e., there is a 
continuous map $u:[0,\infty)\ra \U(A)$ such that 
$\alpha_{t_0}=\lim_{s\rightarrow\infty}\Ad\,u(s)$ and 
 $$
 \lim_{s\rightarrow\infty}\max_{t\in[0,1]}\|\alpha_t(u(s))-u(s)\|=0.
 $$ 
\end{lem}
\begin{pf}
Since $KK(\alpha_{t_0})=KK(\id)$, $\alpha_{t_0}$ is asymptotically 
inner \cite{Phi}, i.e., there is a continuous map 
$v:[0,\infty)\ra\U(A)$ such that 
 $$
 \alpha_{t_0}=\lim_{s\ra\infty}\Ad\,v(s).
 $$               
 
Let $w(s,t)=v(s)\alpha_t(v(s)^*),\ s\in[0,\infty),\ t\in\R$. Then 
for each $s\in[0,\infty)$, the map $t\mapsto w(s,t)$ is an 
$\alpha$-cocycle, i.e., $t\mapsto w(s,t)$ is a continuous function 
into $\U(A)$ such that 
$w(s,t_1+t_2)=w(s,t_1)\alpha_{t_1}(w(s,t_2)),\ t_1,t_2\in\R$. 
Since for each $x\in A$, 
 $$
 \|[w(s,t),x]\|\leq \|\Ad\,v(s)(\alpha_{-t_0-t}(x))-\alpha_{-t}(x)\|
 +\|\Ad\,v(s)(\alpha_{-t_0}(x))-x\|.
 $$
we get, for any $T\gg0$ and  for any $x\in A$, that 
 $$
 \sup_{0\leq t\leq T}\|[w(s,t),x]\|\ra0
 $$
as $s\ra\infty$.   

More specifically let $\F$ be a finite subset of $A$ and $\eps>0$. 
Then there exists an $a>0$ such that if $s\geq a$, then
$\|\Ad\,v(s)(\alpha_{-t_0-t}(x))-\alpha_{-t}(x)\|<\eps/22$ for 
$x\in\F$ and $t\in[0,T]$, which entails  that 
$\|[w(s,t),x]\|<\eps/11$ for $x\in\F$ and $t\in[0,T]$. 

Furthermore, for any bounded interval $I$ of $[0,\infty)$, there 
is a continuous map $z:I\times [0,T]\ra \U(A)$ such that 
 \BE
 z(s,0)&=&1,\\
 z(s,T)&=&w(s,T),\\
 \|z(s,t_1)-z(s,t_2)\|&\leq& (16\pi/3+\eps)|t_1-t_2|/T,\\
 \|[z(s,t),x]\|&<& 10\eps/11, \ 
 \ x\in \F,
 \EE
for $s\in I$ and $t,t_1,t_2\in [0,T]$. (Here we used the estimate 
for a particular construction of $z(s,t)$  that 
 $$\|[z(s,t),x]\|<9\max_{0\leq t_1\leq 
 T}\|[w(s,t_1),x]\|+\eps'$$
for any $\eps'>0$; see \cite{Nak} or 2.7 of \cite{K03}.)  By using 
this $z$,  we get a continuous map $U:I\ra \U(A)$ such that 
 \BE
 \|w(s,t)-U(s)\alpha_t(U(s)^*)\|&\leq & 6\pi|t|/T+\eps,\\
 \|[U(s),x]\|&\leq& \eps,\ \ x\in\F,
 \EE
where we have assumed that $16\pi/3+\eps<6\pi$. 

We recall how $U_T(s)=U(s)$ is defined \cite{K96b}. We define a 
unitary $\tilde{U}_T$ in $C(\R/\Z)\otimes A$ by 
 $$
 \tilde{U}_T(t)=w(s,Tt)\alpha_{T(t-1)}(z_T(s,Tt)^*),
 $$
where $\R/\Z$ is identified with $ [0,1]/\{0,1\}$ and 
$z_T(s,t)=z(s,t)$ is defined above, and we embed $C(\R/\Z)\otimes 
A$ into $A$ {\em approximately} by using the Rohlin property. (If 
$\tau$ is the flow on $C(\R/\Z)$ induced by translations on 
$\R/\Z$, then $\|1\otimes 
w(s,t)-\tilde{U}(\tau_{t/T}\otimes\alpha_{t})(\tilde{U}^*)\|\leq 
6\pi |t| /T$. We find an approximate homomorphism $\phi$ of 
$C(\R/\Z)\otimes A$ into $A$ such that 
$\phi\circ(\tau_{t/T}\otimes \alpha_{t})\approx \alpha_t\phi$ and 
$\phi(1\otimes x) \approx x,\ x\in A$.) Since $z_T$ is defined in 
terms of $w(s,t),\ t\in[0,T]$ and other elements which almost 
commute with them, we may assume that $S\in [0,T]\ra z_S$ is 
continuous; hence that $S\in[0,T]\ra \tilde{U}_S\in 
\U(C(\R/\Z)\otimes A)$ is continuous. Note also that $\tilde{U}_S$ 
commutes with any element to the same degree as $\tilde{U}_T$ does 
with it. Since $\tilde{U}_0=1$, we may thus assume that there is a 
continuous path $(U_t,\ t\in[0,1])$ in the space of continuous 
maps of $I$ into $\U(A)$ such that $U_0(s)=1$, $U_1(s)=U(s)$, and 
$\|[U_t(s),x]\|<\eps,\ x\in\F$. 

Then we set $v_1(s)=U(s)^*v(s)$ for $s\in I$, which satisfies that 
 \BE
 \|\Ad\,v_1(s)(x)-\alpha_{t_0}(x)\|&\leq& 2\eps,\ \ x\in\alpha_{-t_0}(\F), \\   
 \max_{0\leq t\leq1} \|\alpha_t(v_1(s))-v_1(s)\| &\leq &6\pi/T+\eps.
 \EE     

Thus we have shown the following assertion: For any finite subset 
$\F$ of $A$ and $\eps>0$, there exists an $a\in [0,\infty)$ such 
that for any compact interval $I$ of $[a,\infty)$ we find a 
continuous $v_I:I\times [0,1]\ra\U(A)$ such that 
 $$
 \|\Ad\,v_I(s,t)(x)-\alpha_{t_0}(x)\|<\eps,\ \ x\in\F,\ (s,t)\in I\times[0,1],
 $$
and 
 \BE             
 && v_I(s,0)=v(s),\ \ s\in I,\\
 &&\max_{0\leq t\leq1}\|\alpha_t(v_I(s,1))-v_I(s,1)\|<\eps,\ \ s\in 
 I,
 \EE
where $v:[0,\infty)\ra\U(A)$ has been chosen so that 
$\alpha_{t_0}=\lim_{s\rightarrow\infty}\Ad\,v(s)$.

Let $(\F_n)$ be an increasing sequence of finite subsets of $A$ 
such that $\bigcup_n\F_n$ is dense in $A$ and $(\eps_k)$ a 
decreasing sequence in $(0,\infty)$ such that $\lim_k\eps_k=0$. We 
choose an increasing sequence $(a_k)$ in $(0,\infty)$ such that if 
$I$ is a compact interval of $[a_k,\infty)$ then there is a 
continuous $v:I\times [0,1]\ra\U(A)$ such that the above 
conditions are satisfied for $\F=\F_k$ and $\eps=\eps_k$.

Let $a_0=0$ and $I_k=[a_k,a_{k+1}]$ for $k=0,1,2,\ldots$. For each 
$k=1,2,\ldots$ we choose $v_k:I_k\times[0,1]\ra\U(A)$ for $\F_k$ 
and $\eps_k$ as above and define $v_0:I_0\times[0,1]\ra\U(A)$ by 
$v_0(s,t)=v(s),\ s\in I_0$. If $v_{k-1}(a_k)=v_k(a_k)$ for 
$k=1,2,\ldots$, we would be finished by defining a continuous 
function $v:[0,\infty)\ra\U(A)$ with the desired properties in an 
obvious way. But note that $v_k(a_k)v_{k-1}(a_k)^*$ is connected 
to $1$ by a continuous path $(w_k(s),\ s\in[0,1])$ such that 
$w_k(0)=1$, $w_k(1)=v_k(a_k)v_{k-1}(a_k)^*$, and 
 $$
 \|[w_k(s),x]\|<2\eps_{k-1},\ \ x\in \alpha_{t_0}(\F_{k-1})
 $$
for $s\in[0,1]$. By modifying the path $w_k(s),\ s\in[0,1]$, we 
have to impose the condition that $\max_{0\leq 
t\leq1}\|\alpha_t(w_k(s))-w_k(s)\|$ is small; then the path 
$s\mapsto w_k(s)v_{k-1}(s)$ connects $v_{k-1}(a_k)$ with 
$v_k(a_k)$ and has the desired property with respect to $\alpha$. 
Thus it suffices to prove the following lemma by assuming that 
$(\alpha_{t_0}(\F_k))_k$ is sufficiently rapidly increasing and 
$(\eps_k)$ is sufficiently rapidly decreasing. 
\end{pf}        

\begin{lem}\label{A2}
For any finite subset $\F$ of $A$ and $\eps>0$ there exist a 
finite subset $\G$ of $A$ and $\delta>0$ satisfying the following 
condition: If a continuous $v:[0,1]\ra\U(A)$ satisfies that 
 \BE 
 v(0)&=&1,\\
 \|[v(s),x]\|&<&\delta,\ \ s\in[0,1],\ x\in\G,\\
 \|\alpha_t(v(1))-v(1)\|&<&\delta,\ \ t\in[0,1],
 \EE                              
then there exists a continuous $u:[0,1]\ra\U(A)$ such that
 \BE
 u(0)&=&1,\\
 u(1)&=&v(1),\\
 \|[u(s),x]\|&<&\eps,\ \ s\in[0,1],\ x\in\F,\\
 \|\alpha_t(u(s))-u(s)\|&<&\eps,\ \ t\in[0,1],\ s\in[0,1].
 \EE
\end{lem}
\begin{pf} 
Suppose that $v$ satisfies that $v(0)=1$, and
 \BE
 \|[v(s),\alpha_{-t}(x)]\|&<& \delta,\ \ x\in\G,\ t\in[0,1],\\
 \max_{0\leq t\leq1}\|\alpha_t(v(1))-v(1)\|&<&\delta.
 \EE
We define $w(s,t)=v(s)\alpha_t(v(s)^*)$. Then $t\mapsto w(s,t)$ is 
an $\alpha$-cocycle for each $s\in[0,1]$ and satisfies that 
$w(0,t)=1$, and 
 \BE
 \max_{0\leq t\leq 1}\|[w(s,t),x]\|&<&2\delta,\\
 \max_{0\leq t\leq1}\|w(1,t)-1\| &<&\delta.
 \EE  
From the latter condition there are $b,h\in A_{sa}$ such that 
$b\approx0$, $h\approx0$, and 
$w(1,t)=e^{ib}z_t^{(h)}\alpha_t(e^{-ib})$, where $z^{(h)}$ is a 
differentiable $\alpha$-cocycle such that $dz_t/dt|_{t=0}=ih$ 
\cite{K00}. By connecting the $\alpha$-cocycle $t\mapsto w(1,t)$ 
with the trivial $\alpha$-cocycle $1$ by the path of 
$\alpha$-cocycles $s\mapsto (t\mapsto 
e^{isb}z_t^{(sh)}\alpha_t(e^{-isb}))$ and squeezing it around 
$0\in\T=\R/\Z$, we get an $\alpha$-cocycle $W$ in $C(\T)\otimes A$ 
with respect to the flow $\id\otimes\alpha$ such that $W(0,t)=1$ 
and $W(s,t)\approx w(s,t)$. Hence it suffices to show the 
following lemma, because then we find a unitary $Z\in C(\T)\otimes 
A$ with appropriate commutativity such that $Z(0)=1$ and 
$W(\,\cdot\,,t)\approx Z\alpha_t(Z)^*$, and replace $v$ by the 
path  $s\mapsto Z(s)^*v(s)$ which are almost $\alpha$-invariant 
and moves from $v(0)=1$ to $v(1)$. 
\end{pf} 

\begin{lem}
For any finite subset $\F$ of $A$ and $\eps>0$ there exists a 
finite subset $\G$ of $A$ and $\delta>0$ satisfying the following 
condition: Let $\ol{\alpha}=\id\otimes\alpha$ be the flow on 
$C(\T)\otimes A$ and let $t\mapsto W_t$ be an 
$\ol{\alpha}$-cocycle such that $W_t(0)=1$ at 
$0\in\T=[0,1]/\{0,1\}$ and 
 $$
 \max_{0\leq t\leq 1}\|[W_t,1\otimes x]\|<\delta,\ \ x\in\G.
 $$
Then there exists a unitary $Z$ in $C(\T)\otimes A$ such that 
$Z(0)=1$ and
 \BE
 \|[Z,1\otimes x]\|&<&\eps,\ \ x\in\F,\\
 \max_{0\leq t\leq 1}\|W_t-Z\ol{\alpha}_t(Z^*)\|&<&\eps.
 \EE
\end{lem} 
\begin{pf}
We just sketch the proof; see \cite{K96b} or the first part of the 
proof of \ref{A1} for details. 

To meet the last condition we choose $T\in\N$ such that 
$T^{-1}<\eps/6\pi$. Then we impose the condition that $\max_{0\leq 
t\leq 1}\|[W_t,1\otimes x]\|<\delta/T$ for $x\in \bigcup_{-T\leq 
t\leq 0}\alpha_t(\F)$, which can be replaced by a finite subset 
because it is compact. Since $\max_{0\leq t\leq T}\|[W_t,1\otimes 
x]\|<\delta$, we find a continuous path $(U_t,\ t\in[0,T])$ in 
$\U(C(\T)\otimes A)$ such that $U_0=1$, $U_T=W_T$, $U_t(0)=1$, 
$\|U_{t_1}-U_{t_2}\|\leq 6\pi|t_1-t_2|/T$, and 
$\|[U_t,x]\|<9\delta$. By using $W$ and $U$ and the Rohlin 
property for $\alpha$, we define a unitary $Z\in C(\T)\otimes A$ 
such that $Z(0)=1$, $\max_{0\leq t\leq 
1}\|W_t-Z\alpha_t(Z^*)\|<\eps$, and $\|[Z,1\otimes x]\|<10\delta,\ 
x\in \F$. 
\end{pf}      

We also give the following technical results which will be used in 
the next section. We assume that $\alpha$ is a Rohlin flow on $A$ 
as before. 

\begin{lem}\label{B1}
For any finite subset $\F$ of $A$ and $\eps>0$ there exists a 
finite subset $\G$ of $A$ and $\delta>0$ satisfying the following 
condition: If $(u(s),\ s\in[0,1])$ is a continuous path in $\U(A)$ 
such that 
 \BE
 \|[u(s),x]\|&<&\delta,\ \ x\in\G,\\
 \max_{0\leq t\leq 1}\|\alpha_t(u(s))-u(s)\|&<&\delta,\ \ s\in [0,1],
 \EE                 
then there exists a rectifiable path $v(s),\ s\in[0,1]$ such that 
$v(0)=u(0)$, $v(1)=u(1)$, 
 \BE
 \|[v(s),x]\|&<&\eps,\ \ x\in\F,\\
 \max_{0\leq t\leq1}\|\alpha_t(v(s))-v(s)\|&<&\eps,\ \ s\in[0,1],
 \EE
and the length of the path $v$ is less than $17\pi/3$. If 
$\F=\emptyset$, then $\G=\emptyset$ is possible. 
\end{lem} 
\begin{pf} 
Without the conditions with respect to $\alpha$, this is shown in 
\cite{Nak}. 

To define $v$ we use certain elements of $A$ which almost commute 
with $u(s),\ s\in[0,1]$. They are a certain compact subset of 
$\calO_\infty$, which is then embedded centrally in $A$ in 
\cite{Nak}, by using a result due to Kirchberg and Phillips. In 
the present case, to meet the condition of almost 
$\alpha$-invariance, those elements embedded in $A$ should be  
almost invariant under $\alpha$. For this we use the fact that 
$(A'\cap A_\alpha^\omega)^\alpha$ is purely infinite and simple 
\cite{K03}. 

Explicitly we assume that those elements of 
$\calO_\infty=C^*(s_1,s_2,\ldots )$ (before the embedding into 
$A$) is in the linear subspace spanned by a finite number of 
monomials in $s_1,\ldots,s_k$ and their adjoints  for some $k$. We 
find a finite sequence $(T_1,\ldots,T_k)$ of isometries in 
$(A'\cap A_\alpha^\omega)^\alpha$ such that 
$\sum_{i=1}^kT_iT_i^*\lneqq1$. Each $T_i$ is represented by a 
central sequence $(t_i(m))$ of isometries in $A$ such that 
$\sum_{i=1}^kt_i(m)t_i(m)^*\lneqq1$ and 
$\max_{t\in[0,1]}\|\alpha_t(t_i(m))-t_i(m)\|\ra0$ as $m\ra\infty$. 
We then express those elements in terms of $t_1(m),\ldots,t_k(m)$ 
in place of $s_1,\ldots,s_k$ respectively for a sufficiently large 
$m$. Thus we get the required condition involving $\alpha$.
\end{pf}
                                                                 
We will denote by $\delta_\alpha$ the generator of $\alpha$, which 
is a closed derivation from a dense $*$-subalgebra 
$D(\delta_\alpha)$ into $A$. See \cite{BR,Br86,Sak} for the theory 
of generators and derivations. 

\begin{lem}\label{B2}   
Let $(u(s),\ s\in[0,\infty))$ be a continuous path in $\U(A)$ such 
that 
 $$\max_{0\leq t\leq 1}\|\alpha_t(u(s))-u(s)\|
 $$ 
converges to zero as $s\ra\infty$. Then there is a continuous path 
$(v(s),\ s\in[0,\infty))$ of unitaries such that 
$v(s)\in\D(\delta_\alpha)$ and $\delta_\alpha(v(s))$ and 
$u(s)-v(s)$ converge to zero as $s\ra\infty$. 
\end{lem}                   
\begin{pf}
Let $f$ be a non-negative $C^\infty$-function on $\R$ of compact 
support such that the integral is $1$. We set 
 $$
 z(s)=b(s)\int f(b(s)t)\alpha_t(u(s))dt,
 $$
where $b:[0,\infty)\ra(0,\infty)$ is a continuous decreasing 
function such that $\lim_sb(s)=0$, $\|z(s)-u(s)\|<1$, and 
$\|z(s)-u(s)\|\ra0$. Then it follows that $z(s)\in 
\D(\delta_\alpha)$ and 
 $$
 \|\delta_\alpha(z(s))\|\leq b(s)\int|f'(t)|dt, 
 $$
which converges to zero as $s\ra\infty$. We set 
$v(s)=z(s)|z(s)|^{-1}$, which satisfies the required conditions. 
\end{pf}  

\begin{lem} \label{B3}
For a finite subset $\F$ of $A$ and $\eps>0$ there exists a finite 
subset $\G$ of $A$ and $\delta>0$ satisfying the following 
condition. Let $u$ be a unitary in $C[0,1]\otimes A$ such that 
$u(0)=1$, $u(t)\in\D(\delta_\alpha)$, 
$\|\delta_\alpha(u(t))\|<\delta$, and $\|[u(t),x]\|<\delta,\ 
x\in\G$. Then there exist an $h_i\in \D(\delta_\alpha)\cap A_{sa}$ 
for $i=1,2,\ldots,10$ such that 
 \BE
 u(1)&=&e^{ih_1}e^{ih_2}\cdots e^{ih_{10}},\\
 \|h_i\|&<&\pi,\\
 \|\delta_\alpha(h_i)\|&<&\eps,\\
 \|[h_i,x]\|&<&\eps,\ \ x\in\F.
 \EE                                           
If $\F=\emptyset$, then $\G=\emptyset$ is possible. 
\end{lem}  
\begin{pf}
We may assume, by \ref{B1}, that the length of the path $u$ is 
smaller than $17\pi/3<18$. Then we choose $0<s_1<s_2<\cdots 
<s_9<1$ such that $\|u(s_i)-u(s_{i-1})\|<9/5<2$ for 
$i=1,2,\ldots,10$, where $s_0=0$ and $s_{10}=1$. Note that 
 $$
 u(1)=u(s_0)^*u(s_1)\cdot u(s_1)^*u(s_2)\cdots  u(s_9)^*u(s_{10}).
 $$           
Since $\|u(s_{i-1})^*u(s_i)-1\|<9/5$, the spectrum of 
$u(s_{i-1})^*u(s_i)$ is contained in 
 $$
 S=\{e^{i\theta}\ |\ |\theta|<\theta_0\},
 $$                                     
where $\theta_0=\pi-2\cos^{-1}(9/10)<\pi$. Let $\Arg$ denote the 
function $e^{i\theta}\mapsto \theta$ from $S$ onto the interval 
$(-\theta_0,\theta_0)$ and set $h_i=\Arg(u(s_{i-1})^*u(s_i))$. 
Then we have that $\|h_i\|<\pi$ and $u(1)=e^{ih_1}e^{ih_2}\cdots 
e^{ih_{10}}$. We shall show that these $h_i$ satisfy the other 
conditions for a sufficiently small $\delta>0$. 

In general if $v$ is a unitary with $\Spec(v)\subset S$, then 
$h=\Arg(v)$ can be obtained as 
 $$
 h=\frac{1}{2\pi i}\oint_C(\log z) (z-v)^{-1}dz,
 $$
where $\log z$ is the logarithmic function on $\C\setminus 
(-\infty,0]$ with values in $\{z\ |\ |\Im z|<\pi\}$ and $C$ is a 
simple rectifiable path surrounding $S$ in the domain of $\log$. 
We fix $C$ and let $r$ be the distance between $C$ and $S$. Since 
 $$
 \delta_\alpha(h)=\frac{1}{2\pi i}\oint_C\log z(z-v)^{-1}\delta_\alpha(v)
 (z-v)^{-1}dz,
 $$
we have the estimate 
 $$
 \|\delta_\alpha(h)\|\leq (2\pi)^{-1}M|C| r^{-2}\|\delta_\alpha(v)\|,
 $$                                                            
where $M$ is the maximum of $|\log z|,\ z\in C$ and $|C|$ is the 
length of $C$. Similarly we have the estimate $\|[h,x]\|\leq 
(2\pi)^{-1}M|C|r^{-2}\|[v,x]\|$ for any $x\in A$. (See 
\cite{BR,Sak} for details.) Thus we get the conclusion.   
\end{pf}

\section{Proof of Proposition \ref{En}}

We recall that $\calE_n=C^*(s_1,\ldots,s_n)$, where 
$s_1,\ldots,s_n$ are isometries such that 
$e_n^0=1-\sum_{k=1}^ns_ks_k^*$ is a non-zero projection, and that 
$\calJ_n$ is the ideal of $\calE_n$ generated by $e_n^0$. Let 
$\calS=\{1,2,\ldots,n\}^*$ denote the set of all finite sequences 
including an empty sequence, denoted by $\emptyset$. For 
$I=(i_1,i_2,\ldots,i_m)\in \calS$ with $m=|I|$, we set 
$s_I=s_{i_1}s_{i_2}\cdots s_{i_m}$, where $|I|$ is the length of 
$I$; if $|I|=0$ or $I=\emptyset$, then $s_I=1$. It then follows 
that $\{s_Ie_ns_J^*\ |\ I,J\in \calS\}$ forms a family of matrix 
units and spans $\calJ_n$. Thus, in particular, $\calJ_n$ is 
isomorphic to the \cstar\ $\calK$ of compact operators (on an 
infinite-dimensional separable Hilbert space). Hence there is a 
unique (up to unitary equivalence) irreducible representation 
$\pi_0$ of $\calE_n$ such that $\pi_0|\calJ_n$ is non-zero or 
$\pi_0(e_n^0)$ is a one-dimensional projection. We call this 
representation the Fock representation and denote by $\Hil_0$ the 
representation Hilbert space of $\pi_0$.

We recall that the flow $\alpha$ is defined as 
$\alpha_t(s_k)=e^{ip_kt}s_k$ for $k=1,2,\ldots,n$. For 
$I=(i_1,\ldots,i_m)\in\calS$ let $p(I)=\sum_{k=1}^mp_{i_k}\in\R$. 
We set $H_0=\sum_{I\in\calS}p(I)\pi_0(s_Ie_ns_I^*)$, which is a 
well-defined self-adjoint operator on $\Hil_0$. Then it follows 
that $\Ad\,e^{itH_0}\pi_0(x)=\pi_0\alpha_t(x),\ x\in\calE_n$ and, 
by the assumption on $p_1,\ldots,p_n$, that the spectrum of $H_0$ 
is the whole $\R$. Note that 
$\calE_n/\calJ_n\cong\calO_n=C^*(\dot{s}_1,\ldots,\dot{s}_n)$, 
where $\dot{s}_k=s_k+\calJ_n$; we will later on denote $\dot{s}_k$ 
by $s_k$. Note also that $\alpha$ induces a flow on $\calO_n$, 
which we will also denote by $\alpha$. We will denote by $Q$ the 
quotient map of $\calE_n$ onto $\calO_n$.  
                    
\begin{lem}\label{C1}
For any $h\in (\calO_n)_{sa}\cap \D(\delta_\alpha)$ and $\eps>0$ 
there exists a $b\in (\calE_n)_{sa}\cap D(\delta_\alpha)$ such 
that $Q(b)=h$, $\|b\|<\|h\|+\eps$, and 
$\|\delta_\alpha(b)\|<\|\delta_\alpha(h)\|+\eps$. 
\end{lem}
\begin{pf}
Since $Q\alpha_t=\alpha_tQ$ on $\calE_n$, it follows that 
$(1+\delta_\alpha)^{-1}Q=Q(1+\delta_\alpha)^{-1}$, which implies 
that $Q(\D(\delta_\alpha))=\D(\delta_\alpha)$. (We will use the 
same symbol $\delta_\alpha$ for the generator of $\alpha|\calE_n$ 
and of $\alpha|\calO_n$.) 

Thus, for any $h$ as above, there is a $b\in (\calE_n)_{sa}\cap 
\D(\delta_\alpha)$ such that $Q(b)=h$. By $C^\infty$-functional 
calculus we may suppose that $\|b\|<\|h\|+\eps$. 

Since $H_0$ is diagonal, there exists an approximate identity 
$(p_k)$ for $\calJ_n$ consisting of projections in $\calJ_n\cap 
\D(\delta_\alpha)$ such that $\delta_\alpha(p_k)=0$. Since 
$Q(p_k)=0$, we may replace $b$ by $(1-p_k)b(1-p_k)$. We choose a 
$p_k$ such that 
$\|(1-p_k)\delta_\alpha(b)(1-p_k)\|<\|\delta_\alpha(h)\|+\eps$ 
(since $ \|Q(\delta_\alpha(b))\|=\|\delta_\alpha(h)\| 
=\lim_k\|(1-p_k)\delta_\alpha(b)(1-p_k)\|$). Then it follows that 
$b_1=(1-p_k)b(1-p_k)$ belongs to $(\calE_n)_{sa}\cap 
\D(\delta_\alpha)$ and satisfies that $\|b_1\|\leq \|b\|\leq 
\|h\|+\eps$ and 
 $$
 \|\delta_\alpha(b_1)\|= \|(1-p_k)\delta_\alpha(b)(1-p_k)\|
 <\|\delta_\alpha(h)\|+\eps.
 $$ 
Thus $b_1$ satisfies the required conditions.                                                    
\end{pf} 

\begin{lem}\label{C2}
For any finite subset  $\F$ of $\calE_n$ and $\eps>0$ there exists 
a finite subset $\G$ of $\calO_n$ and $\delta>0$ satisfying the 
following condition: If $h\in(\calO_n)_{sa}\cap \D(\delta_\alpha)$ 
such that $\|h\|<\pi$, $\|\delta_\alpha(h)\|<\delta$, and 
$\|[h,x]\|<\delta,\ x\in\G$, then there is a $b\in 
(\calE_n)_{sa}\cap \D(\delta_\alpha)$ such that 
 \BE        
 Q(b)&=&h,\\
 \|b\|&<&\pi+\eps,\\
 \|\delta_\alpha(b)\|&<&\eps,\\
 \|[b,x]\|&<&\eps,\ \ x\in\F,\\
 be_n^0&=&0.
 \EE
\end{lem}
\begin{pf} 
Note that in the above statement we may allow $\F$ and $\G$ to be 
compact subsets instead of finite subsets. 

Note that $\calE_n$ is nuclear as well as $\calO_n$ and $\calJ_n$.
Let $T>0$ be so large that $\pi/T<\eps/2$. Since 
$\F_1=\bigcup_{-T\leq t\leq T}\alpha_t(\F)$ is compact, there is a 
$w=(w_1,\ldots,w_K)\in M_{1K}(\calE_n)$ for some $K\in\N$ such 
that $ww^*=1$ and 
 $$
 \|[wxw^*,a]\|\leq (\eps/\pi)\|x\|, \ \ x\in \calE_n,
 $$                                             
for any $a\in \F_1$, where $wxw^*=\sum_{i=1}^Kw_ixw_i^*$. Note 
that $x\mapsto wxw^*$ is a kind of unital averaging map of $A$ 
into $A$ (this is due to Haagerup \cite{Haa}; see also 
\cite{KOS}).  We set $\G=\{Q(\alpha_t(w_i))\ |\ i=1,2,\ldots,K,\ 
|t|\leq T\}$, which is a compact subset of $\calO_n$. Let 
$\delta\in (0,\eps/2)$ and let $h\in(\calO_n)_{sa}\cap 
\D(\delta_\alpha)$ be such that $\|h\|<\pi$, 
$\|\delta_\alpha(h)\|<\delta$, and $\|[h,x]\|<\delta,\ \ x\in\G$. 
We assume that $\delta>0$ is so small that we get 
 $$
 \|Q(\alpha_t(w))hQ(\alpha_t(w)) -h\|<\eps,\ \ t\in[-T,T].
 $$
Let $b\in(\calE_n)_{sa}\cap \D(\delta_\alpha)$ be such that 
$Q(b)=h$, $\|b\|<\pi$, and $\|\delta_\alpha(b)\|<\delta$. 

We define 
 $$
 b_1=\frac{1}{2T}\int_{-T}^{T}\alpha_t(w)b\alpha_t(w)^*dt.
 $$

Then it follows that $\|b_1\|<\pi$ and $\|Q(b_1)-h\|<\eps$. It 
also follows that $b_1\in (\calE_n)_{sa}\cap\D(\delta_\alpha)$ and 
 $$
 \delta_\alpha(b_1)=(2T)^{-1}(\alpha_T(w)b\alpha_T(w)^*
 -\alpha_{-T}(w)b\alpha_{-T}(w)^*)
 +\frac{1}{2T}\int_{-T}^T\alpha_t(w)\delta_\alpha(b)\alpha_t(w)^*dt,
 $$                                                                 
which implies that $\|\delta_\alpha(b_1)\|<\pi/T+\delta<\eps$. 

Let $a\in\F$. Since 
$\|[\alpha_t(w)b\alpha_t(w)^*,a]\|=\|[w\alpha_{-t}(b)w^*,\alpha_{-t}(a)]\| 
\leq (\eps/\pi)\|b\|<\eps$ for $t\in[-T,T]$, we get that 
$\|[b_1,a]\|<\eps, \ a\in\F$. 

To meet the condition $b_1e_n^0=0$, let $(p_k)$ be an approximate 
identity for $\calJ_n$ in $\calJ_n\cap \D(\delta_\alpha)$ such 
that $p_k\geq e_n^0$, $\delta_\alpha(p_k)=0$ and $\|[p_k,x]\|\ra0$ 
for all $x\in\calE_n$. We replace $b_1$ by $(1-p_k)b_1(1-p_k)$ for 
a sufficiently large $k$.            

In this way we get a $b\in (\calE_n)_{sa}\cap \D(\delta_\alpha)$ 
which satisfies all the required conditions except for $Q(b)=h$; 
instead of which we have that $\|Q(b)-h\|<\eps$. By the previous 
lemma, since $\|\delta_\alpha(Q(b)-h)\|<\eps+\delta$, we get a 
$c\in (\calE_n)_{sa}\cap \D(\delta_\alpha)$ such that 
$Q(c)=h-Q(b)$, $\|c\|<\eps$, and 
$\|\delta_\alpha(c)\|<\eps+\delta$. We may also require that 
$ce_n^0=0$. Thus we can take $b+c$ for $b$, which satisfies the 
required conditions if we start with a smaller $\eps$. 
\end{pf}  
  
Fix $t_0\in\R$. We choose, by \ref{A1} and \ref{B2}, a continuous 
$u:[0,\infty)\ra\U(\calO_n)\cap \D(\delta_\alpha)$ such that 
$\alpha_{t_0}=\lim_{s\rightarrow\infty}\Ad\,u(s)$ and 
$\lim_{s\rightarrow\infty}\delta_\alpha(u(s))=0$. Since the 
unitary group of $\calO_n$ is connected, we may suppose that 
$u(0)=1$.

Let $(\F_k)$ be an increasing sequence of finite subsets of 
$\calE_n$ such that the union $\bigcup_k\F_k$ is dense in 
$\calE_n$ and $(\eps_k)$ a decreasing sequence of positive numbers 
such that $\sum_k\eps_k\equiv\eps\ll 1$. We choose, by \ref{C2}, 
$\G_k=\G$ and $\delta_k=\delta$ for $\F=\F_k$ and $\eps=\eps_k$. 
We may suppose that $(\G_k)$ is increasing and $(\delta_k)$ is 
decreasing to zero. 
                                                   
For the above continuous map $u:[0,\infty)\ra \U(\calO_\infty)\cap 
\D(\delta_\alpha)$, we will choose an increasing sequence $(s_k)$ 
in $[0,\infty)$ with $s_0=0$ such that 
$\|\delta_\alpha(u(s_k)^*u(s_{k+1}))\|$ is sufficiently small for 
$k\geq 0$ and $u(s_{k})^*u(s_{k+1})$ is sufficiently central for 
$k\geq 1$. Specifically, by \ref{B3}, we assume that 
$u(s_k)^*u(s_{k+1})$ has the following factorization: 
 \BE
 u(s_k)^*u(s_{k+1})&=&e^{ih_{k,1}}e^{ih_{k,2}}\cdots e^{ih_{k,10}},\\
 \|h_{ki}\|&<&\pi,\\
 \|\delta_\alpha(h_{ki})\|&<&\delta_k,\\
 \|[h_{ki},x]\|&<& \delta_k,\ \ x\in\G_k,
 \EE                                  
where $\G_0=\emptyset$. Then by \ref{C2} we choose $b_{ki}\in 
(\calE_n)_{sa}\cap \D(\delta_\alpha)$ such that 
 \BE        
 Q(b_{ki})&=&h_{ki},\\
 \|b_{ki}\|&<&\pi+\eps_k,\\
 \|\delta_\alpha(b_{ki})\|&<&\eps_k,\\
 \|[b_{ki},x]\|&<&\eps_k,\ \ x\in\F_k,\\
 \pi_0(b_{ki})\Omega_0&=&0,
 \EE                       
where $\F_0=\emptyset$ and $\Omega_0$ is a unit vector in $\Hil_0$ 
such that $\pi_0(e_n^0)\Omega_0=\Omega_0$. 
 
We set $w_0=1$ and lift $u(s_k)=u(s_{k-1})\cdot 
u(s_{k-1})^*u(s_k)$ for $k\geq 1$ to a unitary in $\calE_n\cap 
D(\delta_\alpha)$ as 
 $$
 w_k=w_{k-1}e^{ib_{k,1}}e^{ib_{k,2}}\cdots e^{ib_{k,10}},
 $$
It then follows that $\Ad\,w_k$ converges on $\calE_n$ as 
$k\ra\infty$ and $Q\circ(\lim\Ad\,w_k)=\alpha_{t_0}\circ Q$. When 
we choose $\G_k$ and $\delta_k$, we should choose them for 
$\F=\F_k\cup \Ad\,w_{k-1}^*(\F_k)$ and $\eps=\eps_k$, which will 
make sure that $\Ad\,w_k^*$ also converges. In this way we have 
$\beta=\lim\Ad\,w_k$ as an automorphism of $\calE_n$, which 
satisfies that 
 $$
 \beta\circ Q=Q\circ\alpha_{t_0}.
 $$
 
Since $\|\alpha_t(w_k)-w_k\|$ is dominated by 
 $$
 \|\alpha_t(w_{k-1})-w_{k-1}\|
 +\sum_{j=1}^{10}\|\alpha_t(e^{ib_{kj}})-e^{ib_{kj}}\|,
 $$                                                       
and since $\|\alpha_t(e^{ib_{kj}})-e^{ib_{kj}}\|\leq 
\|\delta_\alpha(b_{kj})\||t|\leq\eps_k|t|$, we get that 
 $$
 \|\alpha_t(w_k)-w_k\|\leq \|\alpha_t(w_{k-1})-w_{k-1}\|+10\eps_k|t|.
 $$   
Thus if $|t|\leq 1$, then it follows that
 $$
 \|\alpha_t(w_k)-w_k\|\leq 10\sum_{i=1}^k\eps_i<10\eps.
 $$ 
 
In the Fock representation $\pi_0$, since 
$\pi_0(w_k)\Omega_0=\Omega_0$, we have that 
 $$
 \pi_0(w_k)\pi_0(x)\Omega_0=\pi_0(\Ad\,w_k(x))\Omega_0
 $$
converges strongly to $\pi_0(\beta(x))\Omega_0$ for any 
$x\in\calE_n$. Hence $\pi_0(w_k)$ converges strongly to a unitary, 
which we will denote by $W$. Note that 
$\Ad\,W\pi_0(x)=\pi_0\beta(x),\ x\in\calE_n$ and 
$W\Omega_0=\Omega_0$. 

As we have remarked before, the unitary flow $U_t\equiv e^{itH_0}$ 
implements $\alpha$ in $\pi_0$, where $H_0=\sum_{I\in 
\calS}p(I)\pi_0(s_Ie_ns_I^*)$. Since 
$\|U_t\pi_0(w_k)U_t^*-\pi_0(w_k)\|<10\eps$ for $|t|\leq1$, we 
obtain that 
 $$
 \|U_tWU_t^*-W\|\leq 10\eps,\ \ t\in[-1,1].
 $$                                                     
We also have that 
 $$
 \pi_0\beta^{-1}\alpha_t\beta=\Ad(W^*U_tWU_t^*)\pi_0\alpha_t.
 $$ 
 
On the other hand, $\beta^{-1}\alpha_{t}\beta$ is obtained as the 
limit of $\Ad(w_k^*\alpha_t(w_k))\alpha_t$, which implies that
 $$
 \|\alpha_t-\beta^{-1}\alpha_t\beta\|\leq 20\eps,\ \ t\in[-1,1].
 $$
Hence there exists an $\alpha$-cocycle $u$ in $\calE_n$ such that 
$\Ad\,u_t\alpha_t=\beta^{-1}\alpha_t\beta$ and 
$\max_{t\in[-1,1]}\|u_t-1\|$ is at most of order of $400\eps$ (see 
p. 296 of \cite{BR}). Combining the observation in the previous 
paragraph and noting that $\pi_0$ is irreducible,, this implies 
that 
 $$\pi_0(u_t)=c(t)W^*U_tWU_t^*
 $$
for some constant $c(t)\in\T$. Then it follows by simple 
computation that $c(t)=e^{ipt}$ for some $p\in\R$. Thus we know 
that $t\mapsto U_tWU_t^*$ is continuous in norm and that 
$W^*U_tWU_t^*\in \pi_0(\calE_n)$. Since 
$Q\Ad\,u_t\alpha_t=\Ad\,Q(u_t)\alpha_t Q=\alpha_tQ$ on $\calE_n$, 
we also have that $Q(u_t)\in \C1\subset \calO_n$, which implies 
that                                           
 $$
 W^*U_tWU_t^*\in \pi_0(\calJ_n+\C1).
 $$

Since any automorphism of $\calE_n$ is weakly inner in $\pi_0$, 
there is a unitary $V$ on $\Hil_0$ such that 
$\Ad\,V\pi_0=\pi_0\beta^{-1}\alpha_{t_0}$. Since the vector state 
of $\calE_n$ defined through $\Omega_0$ is left invariant under by 
$\beta^{-1}\alpha_{t_0}$, we may define $V$ by 
$V\pi_0(x)\Omega_0=\pi_0(\beta^{-1}\alpha_{t_0}(x))\Omega_0,\ 
x\in\calE_n$. Since $Q\beta^{-1}\alpha_{t_0}=Q$, we have that 
$[V,\pi_0(x)]\in \calK(\Hil_0)$ for $x\in\calE_n$. Regarding $V\in 
M(\calJ_n)$, the multiplier algebra of $\calJ_n$ which identifies 
with $B(\Hil_0)$ through $\pi_0$, and $\calE_n\subset M(\calJ_n)$, 
we have that $\beta^{-1}\alpha_t\beta=\Ad(VU_tV^*U_t^*)\alpha_t$. 
Since $\beta^{-1}\alpha_t\beta=\Ad(W^*U_tWU_t^*)\alpha_t$ from 
above and $VU_tV^*U_t^*\Omega_0=\Omega_0=W^*U_tWU_t^*\Omega_0$, we 
have that 
 $$
 VU_tV^*U_t^*=W^*U_tWU_t^*,\ \ t\in\R,
 $$                            
which implies that $t\mapsto U_tVU_t^*$ is norm-continuous and   
 $$ 
 \|U_tVU_t^*-V\|=\|U_tWU_t^*-W\|.
 $$ 
In particular we have that
 $$
 \max_{t\in [0,1]}\|U_tVU_t^*-V\|\leq 10\eps.
 $$
                                            
Let us denote by $M(\calJ_n)_\alpha$ the C$^*$-subalgebra 
consisting of $x\in M(\calJ_n)$ such that $t\mapsto 
\alpha_t(x)=U_txU_t^*$ is norm-continuous. Summing up the above we 
have shown: 

\begin{lem}\label{C3}                                             
Let $(p_1,\ldots,p_n)$ be a finite sequence in $\R$ and define a 
quasi-free flow $\alpha$ on $\calE_n=C^*(s_1,\ldots,s_n)$ by
 $$
 \alpha_t(s_j)=e^{ip_jt}s_j.
 $$
Suppose that $p_1,\ldots,p_n$ generates $\R$ as a closed 
subsemigroup. (Hence the flow $\dot{\alpha}$ on the quotient 
$\calO_n$ induced by $\alpha$ has the Rohlin property and each 
$\dot{\alpha}_t$ is $\alpha$-invariantly asymptotically inner.) 

Fix $t_0\in\R$. For any $\eps>0$ there exists an automorphism 
$\beta$ of $\calE_n$, a sequence $(w_k)$ in $\U(\calE_n)$, and a 
unitary $V\in M(\calJ_n)_\alpha$ such that $t\mapsto 
V^*\alpha_t(V)\in \calJ_n+\C1$ is an $\alpha$-cocycle, 
 \BE                 
 Q(V)&\in& (\calO_n)',\\
 Q\beta&=& Q\alpha_{t_0},\\ 
 \beta^{-1}\alpha_{t_0}&=&\Ad\,V,\\
 \beta&=&\lim_k\Ad\,w_k,\\
 \max_{|t|\leq1}\|\alpha_t(w_k)-w_k\|&<&\eps,\\
 \max_{|t|\leq 1}\|\alpha_t(V)-V\|&<&\eps,
 \EE                                      
where $Q$ denotes the quotient map of $M(\calJ_n)$ onto 
$M(\calJ_n)/\calJ_n$, which maps $\calE_n$ onto $\calO_n$. 
\end{lem}  
                               
To show that $\alpha_{t_0}$ is $\alpha$-invariantly approximately 
inner, we have to approximate $\beta^{-1}\alpha_{t_0}$ by 
$\Ad\,v$, where $v$ is a unitary in $\calJ_n+1$ which is almost 
$\alpha$-invariant. We will use the following result whose proof 
we will postpone to the next section. 

\begin{lem}\label{C4}
For any $\eps>0$ there exists a $\delta>0$ satisfying the 
following condition: Let $V\in M(\calJ_n)_\alpha$ be a unitary 
such that $Q(V)\in (\calO_n)'$ and 
 $$
 \max_{|t|\leq 1}\|\alpha_t(V)-V\|<\delta.
 $$                                       
Then there exists a rectifiable path $(V_s,\ s\in[0,1])$ in 
$\U(M(\calJ_n)_\alpha)$ such that $V_0=1$, $V_1=V$,  
 \BE
 \|\lambda(Q(V_s))-Q(V_s)\|&<&\eps,\\
 \sup_{s\in[0,1]}\max_{|t|\leq1}\|\alpha_t(V_s)-V_s\|&<&\eps,
 \EE
where $\lambda$ is the unital endomorphism of $M(\calJ_n)/\calJ_n$ 
defined by 
 $$
 \lambda(x)=\sum_{i=1}^nQ(s_i)xQ(s_i)^*.
 $$ 
\end{lem}    

%We define a unital endomorphism $\lambda$ of $M(\calJ_n)\cap 
%\{e_n^0\}'$ into itself by $\lambda(x)=\sum_is_ixs_i^*+xe_n^0$. 
%Note that $Q\lambda=\lambda Q$, where the latter $\lambda$ is 
%defined in the previous lemma. 

The following is a key lemma for the proof of Proposition 
\ref{En}. 

\begin{lem} \label{C5}
For any $\eps>0$ there exists a $\delta>0$ satisfying the 
following condition: If $V$ is a unitary in $M(\calJ_n)_\alpha$ 
such that $\beta=\Ad\,V$ is an automorphism of $\calE_n$, 
$Q(V)\in(\calO_n)'$, and 
$\max_{|t|\leq1}\|\alpha_t(V)-V\|<\delta$, then there is a unitary 
$v$ in $\calJ_n+1$ such that 
 \BE
 \|\beta(s_i)-vs_iv^*\|&<&\eps,\ \ i=1,2,\ldots,n,\\
  \max_{|t|\leq 1}\|\alpha_t(v)-v\|&<&\eps.
 \EE 
\end{lem}
\begin{pf}
We define a non-unital endomorphism $\lambda$ of 
$M(\calJ_n)_\alpha$ by $\lambda(x)=\sum_{i=1}^ns_ixs_i^*$. Note 
that $Q\lambda=\lambda Q$, where the latter $\lambda$ is the 
unital endomorphism defined in \ref{C4}. 

Let $\eps>0$. We choose $\delta>0$ so small that we find a 
continuous path $(V_s,\ s\in[0,1])$ in $M(\calJ_n)_\alpha$ such 
that $V_0=V$, $V_1=1$, $\|\lambda(V_s)-V_s+\calJ_n\|<\eps$, and 
$\max_{|t|\leq 1}\|\alpha_t(V_s)-V_s\|<\eps$. 

Let $\eps'>0$. We choose an increasing sequence $(\mu_k)$ in 
$[0,1)$ such that $\mu_0=0$, $\lim_k \mu_k=1$, and 
$\|V_{\mu_k}-V_{\mu_{k+1}}\|<\eps'$. We will denote $V_{\mu_k}$ by 
$V_k$ below. 

Note that $p_m=\sum_{|I|\leq m}s_Ie_n^0s_I^*$ is an 
$\alpha$-invariant projection in $\calJ_n$ for each $m\in\N$ and 
that $(p_m)$ forms an approximate identity for $\calJ_n$. 

Since $V_ss_kV_s^*=V_s\lambda(V_s^*)s_k$ and 
$\|1-V_s\lambda(V_s^*)+\calJ_n\|<\eps$, we have that 
$\|V_ss_kV_s^*-s_k+\calJ_n\|<\eps$. Since $\{V_ss_kV_s^*-s_k\ |\ 
s\in[0,1]\}$ is compact, we have a projection $p\in\calJ_n$ such 
that $\|(V_ss_kV_s^*-s_k)(1-p)\|<\eps$ and 
$\|(1-p)(V_ss_kV_s^*-s_k)\|<\eps$ for $s\in[0,1]$ and 
$k=1,\ldots,n$. We may suppose that $p=p_m$ for some $m$. From the 
convex combinations of $(p_m)$, we find an approximate unit 
$(e_k)$ in $\calJ_n$ such that $e_0=0\leq p\leq e_1$, 
$\alpha_t(e_k)=e_k$, $e_{k+1}e_k=e_k$, and 
 \BE
  \|[(e_{k+1}-e_k)^{1/2},s_i]\|&<&\eps'2^{-k-1},\ \ i=1,2,\ldots,n,\\
  \|[(e_{k+1}-e_k)^{1/2},V_{j}]\|&<&\eps',\ \ j\leq k+1,\\ 
  \|[(e_k-e_k^2)^{1/2},V_{j}]\|&<&\eps',\ \ j\leq k+1.
 \EE                   
Since $-1\leq 
\sum_{k=0}^K(e_{k+1}-e_k)^{1/2}x_k(e_{k+1}-e_k)^{1/2}\leq1$ for 
any $K$ and any $x_k=x_k^*$ with $\|x_k\|\leq1$, we can define  
 $$
 z=\sum_{k=0}^\infty (e_{k+1}-e_k)^{1/2}V_{k}(e_{k+1}-e_k)^{1/2},
 $$                                                                
which converges in the strict topology in $M(\calJ_n)$ and has 
$\|z\|\leq2$.  Since $\|V_k-1\|\ra0$, it follows that 
$z-1\in\calJ_n$. We claim that $z$ is close to a unitary by 
writing 
 \BE
 zz^*&=&\sum_k(e_{k+1}-e_k)^{1/2}V_{k}(e_{k+1}-e_k)V_k^*(e_{k+1}-e_k)^{1/2}\\
  &+&\sum_k(e_{k+1}-e_k)^{1/2}V_k(e_{k+1}-e_{k+1}^2)^{1/2}V_{k+1}^*(e_{k+2}-e_{k+1})^{1/2}\\
  &+&\sum_k(e_{k+2}-e_{k+1})^{1/2} V_{k+1}(e_{k+1}-e_{k+1}^2)^{1/2}V_k^*(e_{k+1}-e_k)^{1/2},                                                                                    
  \EE    
where we have used that  
$(e_{k+1}-e_k)^{1/2}(e_{j+1}-e_j)^{1/2}=0$ if $|k-j|>1$ and 
$(e_{k+1}-e_k)^{1/2}(e_{k+2}-e_{k+1})^{1/2}=(e_{k+1}-e_{k+1}^2)^{1/2}$.  
By splitting out each summation into the sum over even integers 
and the sum over odd integers and noting that 
 $$\|V_{k}(e_{k+1}-e_k)V_{k}^*-(e_{k+1}-e_k)\|<2\eps'
 $$ 
and 
 $$\|V_k(e_{k+1}-e_{k+1}^2)^{1/2}V_{k+1}^*-(e_{k+1}-e_{k+1}^2)^{1/2}\|<2\eps',
 $$ 
we estimate  
 $$
 \|zz^*-1\|<2(2\eps'+2\eps'+2\eps')=12\eps'.
 $$                                  
Similarly we get that $\|zz^*-1\|<12\eps'$. Thus if $\eps'$ is 
sufficiently small, then $v=z|z^*z|^{-1/2}$ is a unitary in 
$\calJ_n+1$ and satisfies that 
 $$\|v-z\|\leq \|v\|\|1-|z^*z|^{1/2}\|<\eps'',
 $$
with $\eps''\approx 6\eps'$. 
                                                         
Since  
 $$
 \|s_iz^*- 
 \sum_{k=0}^\infty (e_{k+1}-e_k)^{1/2}s_iV_k^*(e_{k+1}-e_k)^{1/2}\|<\eps',
 $$                                                         
we have that 
 \BE
 zs_iz^*&\approx& \sum_k (e_{k+1}-e_k)^{1/2}(e_{k+1}-e_k) V_k
 s_iV_k^*(e_{k+1}-e_k)^{1/2}\\
 &+& \sum_k(e_{k+1}-e_k)^{1/2}(e_{k+1}-e_{k+1}^2)^{1/2}V_{k+1}
 s_iV_{k+1}^*(e_{k+2}-e_{k+1})^{1/2}\\
 &+&\sum_k(e_{k+2}-e_{k+1})^{1/2}(e_{k+1}-e_{k+1}^2)^{1/2}V_{k}
 s_iV_k^*(e_{k+1}-e_k)^{1/2},
 \EE                                  
with the norm difference less than $\|z\|\eps'+12\eps'\leq 
14\eps'$. By using $\|(V_ks_iV_k^*-s_i)(1-p)\|<\eps$ etc. for all 
terms except for the $k=0$ term of the first summation, we have 
that 
 \BE
 zs_iz^*&\approx& e_1^{3/2}Vs_iV^*e_1^{1/2}
 +\sum_{k\geq1}(e_{k+1}-e_k)^{3/2}s_i(e_{k+1}-e_k)^{1/2}\\
 &+& \sum_{k\geq 0}(e_{k+1}-e_k)^{1/2}(e_{k+1}-e_{k+1}^2)^{1/2}
 s_i(e_{k+2}-e_{k+1})^{1/2}\\
 &+&\sum_{k\geq 0}(e_{k+2}-e_{k+1})^{1/2}(e_{k+1}-e_{k+1}^2)^{1/2}
  s_i(e_{k+1}-e_k)^{1/2}
 \EE
with the norm difference less than $14\eps'+6\eps$. Since 
$\|[Vs_iV^*,e_1^{1/2}]\|<3\eps'$ and 
$\|[(e_{k+1}-e_k)^{1/2},s_i]\|<\eps'2^{-k-1}$, we have that 
 \BE
 zs_iz^* &\approx&
 e_1^2Vs_iV^*+(1-e_1^2)s_i, 
 \EE                       
with the norm difference less than 
$14\eps'+6\eps+3\eps'+3\eps'=6\eps+20\eps'$. Since 
$e_1^2Vs_iV^*+(1-e_1^2)s_i-\beta(s_i)=(1-e_1^2)(s_i-\beta(s_i))$, 
it follows that $\|zs_iz^*-\beta(s_i)\|<7\eps+20\eps'$, which 
implies that 
$\|vs_iv^*-\beta(s_i)\|<\eps''(\|z\|+1)+7\eps+20\eps'\approx 
7\eps+32\eps'$. 

On the other hand we have that $\|\alpha_t(z)-z\|\leq 
2\sup_k\|\alpha_t(V_k)-V_k\|$, which implies that 
$\max_{|t|\leq1}\|\alpha_t(v)-v\|<2\eps+2\eps''$. This concludes 
the proof.
\end{pf}    

\medskip
\noindent {Proof of Proposition \ref{En}}

This follows by combining \ref{C3}, \ref{C4} and \ref{C5}. 

Let $t_0\in\R$ and $\eps>0$. Then we choose $\delta>0$ for $\eps$ 
as in \ref{C5}. By \ref{C3} we have an automorphism $\beta$ of 
$\calE_n$ and a unitary $V\in M(\calJ_n)_\alpha$ such that 
$Q\beta=\alpha_{t_0}Q$ and $\beta^{-1}\alpha_{t_0}=\Ad\,V$, 
$\max_{|t|\leq1}\|\alpha_t(V)-V\|<\delta$, etc. Then by \ref{C5} 
we get a unitary $v\in\calJ_n+1$ such that 
$\|\beta^{-1}\alpha_{t_0}(s_i)-vs_i v^*\|<\eps$ and $\max_{|t|\leq 
1}\|\alpha_t(v)-v\|<\eps$. Since $\beta=\lim_k\Ad\,w_k$, we get 
that $\|\alpha_{t_0}(s_i)-w_kvs_iv^*w_k\|<\eps$ for a large $k$ 
and $\max_{|t|\leq1} \|\alpha_t(w_kv)-w_kv\|<\eps+\delta$ from the 
properties imposed on $(w_k)$. Thus we can conclude that 
$\alpha_{t_0}$ is $\alpha$-invariantly approximately inner.

\section{Proof of Lemma \ref{C4}}

We recall that $H_0$ is a self-adjoint operator on $\Hil_0$ with 
$\Spec(H_0)=\R$ and that we let $\alpha_t=\Ad\,e^{itH_0}$ on 
$M(\calJ_n)_\alpha$, where $M(\calJ_n)$ is identified with 
$B(\Hil_0)$ and $M(\calJ_n)_\alpha$ is the largest 
C$^*$-subalgebra of $M(\calJ_n)$ on which $t\mapsto 
\Ad\,e^{itH_0}$ is strongly continuous. The following can be 
proved by adopting the arguments in \cite{BEEK,K01}. 

\begin{lem}  \label{U1}
For any $\eps>0$ there exists a $\delta>0$ satisfying the 
following condition: If $Z\in \U(M(\calJ_n)_\alpha)$ satisfies 
that $\max_{|t|\leq 1}\|\alpha_t(Z)-Z\|<\delta$, then there exists 
a rectifiable  path $(Z_s,\ s\in[0,1])$ in $\U(M(\calJ_n)_\alpha)$ 
such that $Z_0=1$, $Z_1=Z$, 
 $$
      \sup_{s\in [0,1]}\max_{|t|\leq1}\|\alpha_t(Z_s)-Z_s\|<\eps, 
 $$
and the length of $(Z_s,\ s\in[0,1])$ is less than $2\pi+\eps$. 
\end{lem} 
\begin{pf}
We choose a non-negative $C^\infty$-function on $\R$ such that 
$\hat{f}(0)=1$ and $\supp\hat{f}\subset (-\delta_0,\delta_0)$ for 
a small $\delta_0>0$ and modify $Z$ by 
 $$
 X=\int f(t)\alpha_t(Z)dt\in M(\calJ_n)_\alpha
 $$                      
which is still close to $Z$, e.g., $\|X-Z\|<\mu$, where $\mu$ can 
be arbitrarily close to zero depending on $\delta$. Let $E$ be the 
spectral measure of $H_0$ and define, for any $k\in\Z$,
 $$
 E_k=E[2k\delta_0,2(k+1)\delta_0)\in M(\calJ_n)_\alpha,
 $$                               
which is a projection of infinite rank. Noting that the 
$\alpha$-spectrum of $X$ is contained in $(-\delta_0,\delta_0)$, 
we have, as in the proof of 4.1 of \cite{K01}, that 
 \BE
 XE_kX^*&\leq & E[(2k-1)\delta_0,(2k+3)\delta_0),\\ 
 XE_kX^*-(XE_kX^*)^2&\leq & 2\mu\cdot XE_kX^*\leq 2\mu1,\\
 \Spec(XE_kX^*)&\subset & \{0\}\cup [1-2\mu,1],\\                       
 Y_k^+-(Y_k^+)^2&\leq & 6\mu1,\\
 \Spec(Y_k^+)&\subset& [0,6\mu']\cup[1-6\mu',1],\\
 Y_k^- -(Y_k^-)^2&\leq & 6\mu,\\
 \Spec(Y_k^-)&\subset& [0,6\mu']\cup[1-6\mu',1],
 \EE
where $\mu\approx\mu'$ and
 \BE
 Y_k^+&=&E[(2k+1)\delta_0,(2k+3)\delta_0)XE_kX^*E[(2k+1)\delta_0,(2k+3)\delta_0),\\
 Y_k^-&=&E[(2k-1)\delta_0,(2k+1)\delta_0)XE_kX^*E[(2k-1)\delta_0,(2k+1)\delta_0).
 \EE
We define $F_k^{\pm}$ to be the spectral projection of $Y_k^{\pm}$ 
corresponding to $[1-6\mu',1]$ and let $F_k$ denote the projection 
$F_k^++F_k^-$. We note that these projections are of infinite 
rank, because $E(I)$ is a projection of infinite rank for any 
non-empty open subset $I$ of $\R$, and that 
 \BE
 F_k^+&\leq& E[(2k+1)\delta_0,(2k+3)\delta_0),\\
 F_k^-&\leq &E[(2k-1)\delta_0,(2k+1)\delta_0),\\
 \|F_k-XE_kX^*\| &<&14\mu.
 \EE
We define 
 \BE
 G_k^+&=& E[(2k+1)\delta_0,(2k+3)\delta_0)-F_{k+1}^-,\\
 G_k^-&=& E[(2k-1)\delta_0,(2k+1)\delta_0)-F_{k-1}^+,\\
 G_k&=&G_k^-+G_k^+.
 \EE
We note that $G_k^{\pm}$ and $G_k$ are projections of infinite 
rank and that
 $$
 \|XE_kX^*-G_k\|<34\mu.
 $$   
By using the facts that $\sum_kE_k=1$ and 
 $$
 \sum_kF_{2k}+\sum_kG_{2k-1}=1,
 $$                            
we check that
 $$ 
 W=\sum_k F_{2k}XE_{2k}+\sum_{k}G_{2k-1}XE_{2k-1},
 $$
is close to a unitary; we denote by $V$ the unitary part of the 
polar decomposition of $W$. Since 
$\|\alpha_t(F_{2k}XE_{2k})-F_{2k}XE_{2k}\|\leq 
\|\alpha_t(X)-X\|+4\delta_0|t|$,  
$\|\alpha_t(G_{2k-1}XE_{2k-1})-G_{2k-1}XE_{2k-1}\|\leq 
\|\alpha_t(X)-X\|+4\delta_0|t|$, and 
 $$
 \|\alpha_t(W)-W\|\leq 
 \sup_k\|\alpha_t(F_{2k}XE_{2k})-F_{2k}XE_{2k}\|+
 \sup_k\|\alpha_t(G_{2k-1}XE_{2k-1})-G_{2k-1}XE_{2k-1}\|,
 $$
we have that $W,V\in M(\calJ_n)_\alpha$. We estimate that 
$\|Z-V\|<10\mu''$, where $\mu''\approx\mu^{1/2}$. We should also 
note that 
 \BE
 VE_{2k}V^*&=& F_{2k}\leq E[(4k-1)\delta_0,(4k+3)\delta_0),\\
 VE_{2k-1}V^*&=& G_{2k-1}\leq E[(4k-3)\delta_0,(4k+1)\delta_0).
 \EE
It then follows that the $\alpha$-spectrum of $V$ is contained in 
$[-3\delta_0,3\delta_0]$.

Since $\|VZ^*-1\|<10\mu''$, there is a self-adjoint $B\in 
M(\calJ_n)_\alpha$ such that $VZ^*=e^{iB}$ and $\|B\|$ is of the 
order $10\mu''$. Then the path $\gamma:[0,1]\ni s\mapsto e^{isB}Z$ 
goes from $Z$ to $V$ of length $\|B\|$. Since 
$\|\alpha_t(e^{isB}Z)-e^{isB}Z\|\leq 
s\|\alpha_t(B)-B\|+\|\alpha_t(Z)-Z\|$, we have that 
$\max_s\|\alpha_t(\gamma_s)-\gamma_s\|\ra0$ as $t\ra0$ and 
 $$
 \sup_s\max_{|t|\leq1}\|\alpha_t(\gamma_s)-\gamma_s\|
 \leq 2\|B\|+ \max_{|t|\leq1}\|\alpha_t(Z)-Z\|,
 $$                                            
which is arbitrarily small. 

Then we find a rectifiable path $(\zeta_t,\ t\in[0,1])$ in 
$\U(M(\calJ_n))$ such that $\zeta_t$ commutes with 
$E[(2k-1)\delta_0,(2k+1)\delta_0)$ for all $k\in\Z$ and 
 \BE
 \Ad\zeta_1(F_{2k}^+)&=&E[(4k+1)\delta_0,(4k+2)\delta_0),\\
 \Ad\zeta_1(F_{2k}^-)&=&E[4k\delta_0,(4k+1)\delta_0),
 \EE
and the length of $\zeta$ is at most $\pi$. Then we get that
 $$
 \Ad(\zeta_1V)E_{2k}=E_{2k}.
 $$                         
Since $\Ad\zeta_1(G_{2k+1}^-)=E[(4k+2)\delta_0,(4k+3)\delta_0)$ 
etc., we also get that
 $$
 \Ad(\zeta_1V)E_{2k+1}=E_{2k+1}.
 $$ 
Note that the $\alpha$-spectrum of $\zeta_t$ is contained in 
$[-2\delta_0,2\delta_0]$, which implies that $\zeta_t\in 
M(\calJ_n)_\alpha$ and $\max_{|t|\leq 
1}\|\alpha_t(\zeta_s)-\zeta_s\|\leq 2\delta_0$. 
                               
There is a rectifiable path $(\eta_t,\ t\in[0,1])$ in 
$\U(M(\calJ_n))$ such that $\eta_t$ commutes with 
$E[(2k-1)\delta_0,(2k+1)\delta_0)$ for all $k\in\Z$, $\eta_0=1$, 
$\eta_1=(\zeta_1V)^*$, and the length of $\eta$ is at most $\pi$. 
Note also that the $\alpha$-spectrum of $\eta_t$ is contained in 
$[-2\delta_0,2\delta_0]$, which implies that $\eta_t\in 
M(\calJ_n)_\alpha$ and 
$\max_{|t|\leq1}\|\alpha_t(\eta_s)-\eta_s\|\leq 2\delta_0$.  
  
Combining these paths we get the conclusion.
\end{pf}          

Next we will prove another version of the above lemma; we will 
replace $M(\calJ_n)_\alpha$ by the C$^*$-tensor product 
$C(\T)\otimes M(\calJ_n)_\alpha$ with the flow $\id\otimes 
\alpha$, which will sometimes be denoted by $\alpha$. 
                                           
For the preparation we present the following two lemmas, which  
are just concerned with $B(\Hil)$, the bounded operators on an 
infinite-dimensional Hilbert space $\Hil$, without a flow on it. 
We note that the C$^*$-tensor product $C(\T)\otimes B(\Hil)$ 
identifies with the norm-continuous functions on $\T$ into 
$B(\Hil)$. 

\begin{lem}    \label{U2}
Let $E,F,P$ be projections in $C(\T)\otimes B(\Hil)$ such that 
$E(s),F(s),P(s)$ are of infinite rank, $P(s)=P(0)$ for 
$s\in\T=[0,1]/\{0,1\}$, and $EP=0=FP$. Then there is a rectifiable 
path $(U_t,\ t\in[0,1])$ in $\U(C(\T)\otimes B(\Hil))$ such that 
$U_0=1$, $\Ad\,U_1(E)=F$, and the length of $U$ is at most $\pi$. 
If furthermore $E(0)=F(0)$ at $0\in\T$, then the condition 
$U_t(0)=1,\ t\in[0,1]$ can be imposed. 

Moreover if $E_k,F_k,P_k$ are such a triple of projections in 
$C(\T)\otimes B(\Hil)$ for each $k\in\N$ such that $\T\ni s\mapsto 
\bigoplus_kE_k(s)$ and $\T\ni s\mapsto \bigoplus _k F_k(s)$ are 
continuous in $\Pi_{k=1}^\infty C(\T)\otimes B(\Hil)$, there are 
rectifiable paths $(U_t^k,\ t\in [0,1])$ in $\U(C(\T)\otimes 
B(\Hil))$ as above such that $\T\ni s\mapsto \bigoplus _k 
U_t^k(s)$ is continuous and the length of $[0,1]\ni t\mapsto 
\bigoplus_kU^k_t$ is at most $\pi$.  
\end{lem} 
\begin{pf}                                
We have expressed the base space $\T$ as $[0,1]/\{0,1\}$. Since 
$t\mapsto E(t)$ is norm-continuous, there exists a norm-continuous 
$V:[0,1]\mapsto \U(B(\Hil))$ such that $V(0)=1$ and 
 $$\Ad\, V(t)(E(0))=E(t), \ t\in[0,1].
 $$
More specifically we choose a finite number of points 
$t_0=1<t_1<t_2<\cdots <t_{m-1}<t_m=1$ such that 
$\|E(s)-E(t)\|<1/2$ for $s,t\in[t_{i-1},t_{i}]$ for 
$i=1,2,\ldots,m$. We define, for $t\in [t_{i-1},t_i]$, 
 $$
 Z_t^i=E(t)E(t_{i-1})+(1-E(t))(1-E(t_{i-1})).
 $$
Since $\|(Z_t^i)^*Z_t^i-1\|<1/2$ etc. the polar decomposition of 
$Z_t^i$ gives a unitary $V_t^i=Z_t^i|Z_t^i|^{-1}$. Then we define 
for $t\in [t_{i-1},t_i]$ 
 $$                     
 V_t=V_t^iV_{t_{i-1}}^{i-1}V_{t_{i-2}}^{i-2}\cdots V_{t_1}^1.
 $$                                        
Since $\Ad\,V^i_t(E(t_{i-1}))=E(t)$, we get that 
$\Ad\,V(t)(E(0))=E(t)$. It is obvious that $t\mapsto V_t$ is 
norm-continuous. But more is true.

Note that $\|V_t^i-V_s^i\|\leq 
2^{1/2}\|E(t)-E(s)\|+\||Z_t^i|^{-1}-|Z_s^i|^{-1}\|$ for $s,t\in 
[t_{i-1},t_i]$. 

Let $\calS$ be the set of self-adjoint elements $h\in B(\Hil)$ 
such that $\Spec(h)\subset [1/2,1]$; then 
$|Z_t^i|^2=(Z_t^i)^*Z_t^i\in \calS$. The map $\calS\ni h\mapsto 
h^{-1/2}$ is uniformly continuous in $h\in\calS$. Since 
$\|\|Z_t^i|^2-|Z_s^i|^2|\leq \|E(t)-E(s)\|$, the continuity of 
$t\mapsto V_t$ only depends on the continuity of $t\mapsto E(t)$ 
on each $[t_{i-1},t_i]$, i.e., there is a non-decreasing 
continuous function $\varphi_1$ on $[0,1]$ such that 
$\varphi_1(0)=0$ and $\|V_s-V_t\|\leq \varphi_1(\|E(s)-E(t)\|)$ 
for $s,t\in [t_{i-1}-t_i]$ with all $i$. Since $[0,1]\ni t\mapsto 
E(t)$ is uniformly continuous, there is a non-decreasing 
continuous function $\varphi_2$ on $[0,1]$ into $[0,1]$ such that 
$\varphi_2(0)=0$ and $\|E(s)-E(t)\|\leq \varphi_2(|s-t|)$. Let 
$\varphi=\varphi_1\circ\varphi_2$. Combining these estimates, we 
have that $\|V_s-V_t\|\leq \varphi(|s-t|)$ for $s,t\in 
[t_{i-1},t_i]$. 

Let $\Delta=\min_{1\leq i\leq m}\{|t_i-t_{i-1}|\}$ and assume, by 
replacing $\varphi$ by a bigger one if necessary, that 
$\varphi(s)+\varphi(t)\leq \varphi(s+t)$ for $s,t>0$ with $s+t\leq 
1$. If $s,t\in [0,1]$ satisfies that $0<s-t<2\Delta$, then either 
$t_{i-1}<t<t_i<s<t_{i+1}$ or $t_{i-1}\leq t<s\leq t_i$ for some 
$i$. In the former case $\|V_s-V_t\|\leq 
\|V_s-V_{t_i}\|+\|V_{t_i}-V_t\|\leq 
\varphi_1(\|E(s)-E(t_i)\|)+\varphi_1(\|E(t_i)-E(t)\|)\leq 
\varphi(s-t_i)+\varphi(t_i-t)\leq \varphi(s-t)$. With the same 
estimate for the latter case we have that $\|V_s-V_t\|\leq 
\varphi(|s-t|)$ if $|s-t|<2\Delta$. By redefining $\varphi$ on 
$[2\Delta,1]$, we may suppose that $\|V_s-V_t\|\leq 
\varphi(|s-t|)$ is valid for all $s,t\in[0,1]$. 

We define a function $F$ on $\T=\{z\in \C\ |\ |z|=1\}$ by 
$F(e^{it})=t$ for $t\in (-\pi,\pi]$ and let $H=F(V(1))$. Since 
$[V(1),E(0)]=0$, we have that $[H,E(0)]=0$. We replace $V$ by $V': 
t\mapsto V_te^{-itH}$. Then we get that $V'(1)=1=V'(0)$, 
$\Ad\,V'(t)(E(0))=E(t)$, and the continuity of $t\mapsto V_t'$ 
depends only on the continuity of $t\mapsto E(t)$. If we replace 
$\varphi$ by $s\mapsto \varphi(s)+\pi s$, then we have that 
$\|V'(s)-V'(t)\|\leq \varphi(|s-t|)$.                                         

In this way we construct $V,W\in \U(C(\T)\otimes 
B((1-P(0))\Hil))$, which satisfies that $V(1\otimes E(0))V^*=E$, 
$V(0)=1$, $W(1\otimes F(0))W^*=F$, and $W(0)=1$. Moreover we may 
assume that there is a non-decreasing continuous function 
$\varphi$ on $[0,1]$ such that $\varphi(0)=0$ and
 $$
 \|V(s)-V(t)\|\leq \varphi(|s-t|),\ \ \|W(s)-W(t)\|\leq 
 \varphi(|s-t|).
 $$
Note that $\varphi$ depends only on the choice of $t_i$'s and the 
continuity of $s\mapsto E(s)$ and of $s\mapsto F(s)$.
 
If $E(0)=F(0)$, the unitary $Z\equiv WV^*$ in $C(\T)\otimes 
B((1-P(0))\Hil)$ satisfies that $\Ad(Z)(E)=F$ and $Z(0)=1-P(0)$. 
Let $Y\in B(\Hil)$ be such that $Y^*Y=P(0)$ and $YY^*=1-P(0)$ and 
let 
 $$
 X_t=\cos(\pi t/2)+(Y-Y^*)\sin(\pi t/2),\ \ t\in[0,1].
 $$
Then $X_t$ is a unitary on $\Hil$. The path $U:t\mapsto 
(Z+P)(1\otimes X_t)(Z^*+P)(1\otimes X_t^*)$ in $\U(C(\T)\otimes 
B(\Hil))$ satisfies that $U_0=1$, $U_1=Z+Y^*Z^*Y$, $U_t(0)=1$, and 
 $$
 \|U_{t_1}-U_{t_2}\|\leq 2\|X_{t_1}-X_{t_2}\|\leq \pi|t_1-t_2|.
 $$                      
Moreover $U_t$ satisfies that
 $$
 \|U_t(s_1)-U_t(s_2)\|\leq 2\varphi(|s_1-s_2|).
 $$
Since $\Ad\,U_1(E)=F$, this concludes the proof for the case 
$E(0)=F(0)$. 

In the case $E(0)\not=F(0)$ we choose a unitary $T$ on 
$(1-P(0))\Hil$ such that $TE(0)T^*=F(0)$. Then we can proceed as 
above with $Z=W(1\otimes T)V^*$.

Suppose that $E_k,F_k,P_k$ are given as in the statement. Let 
$E_\infty=\bigoplus_kE_k$ and $F_\infty=\bigoplus_kF_k$ in 
$\Pi_{k=1}^\infty C(\T)\otimes B(\Hil)$. Then we choose a finite 
number of points $t_0=0<t_1<t_2<\cdots <t_m=1$ such that 
$\|E_\infty(s)-E_\infty(t)\|<1/2$ and 
$\|F_\infty(s)-F_\infty(t)\|<1/2$ for $s,t\in [t_{i-1}-t_i]$ with 
all $i$. By using these points in $[0,1]$ we construct a path 
$U^k:[0,1]\ra\U(C(\T)\otimes B(\Hil))$ as above. Then there is a 
non-decreasing continuous function $\varphi$ on $[0,1]$ such that 
$\varphi(0)=0$ and $\|U^k_t(s_1)-U^k_t(s_2)\|\leq 
2\varphi(|s_1-s_2|)$. Note that we can define $\varphi$ based on 
the functions $s\mapsto E_\infty(s)$ and $s\mapsto F_\infty(s)$, 
i.e., $\varphi$ is independent of $k$.  Hence it follows that 
$\T\ni s\mapsto U^\infty_t(s)\equiv \bigoplus_kU^k_t(s)$ is 
continuous for each $t\in[0,1]$. Since 
$\|U_{t_1}^k-U_{t_2}^k\|\leq \pi|t_1-t_2|$, we also get that 
$\|U^\infty_{t_1}-U^\infty_{t_2}\|\leq \pi |t_1-t_2|$. This 
concludes the proof. 
\end{pf}
                     
\begin{lem} \label{U3}
Let $U\in \U(C(\T)\otimes B(\Hil))$. If there is a projection 
$E\in B(\Hil)$ such that $E$ and $1-E$ are of infinite rank and 
$[U,1\otimes E]=0$, then there is a rectifiable path $(V_t,\ 
t\in[0,1])$ in $\U(C(\T)\otimes B(\Hil))$ such that $V_0=1$, 
$V_1=U$, $\|V_t(s_1)-V_t(s_2)\|\leq 6\|U(s_1)-U(s_2)\|$, and the 
length of $V$ is at most $3\pi$. 
\end{lem}
\begin{pf}                                                        
Let $Y\in B(\Hil)$ be a partial isometry such that $Y^*Y=E$ and 
$YY^*=1-E$ and let $X_t=\cos(\pi t/2)+(Y-Y^*)\sin(\pi t/2)$ as in 
the proof of the previous lemma. Then the map $[0,1]\ni t\mapsto 
(U(1\otimes E)+1\otimes (1-E))(1\otimes X_t)(U(1\otimes 
(1-E))+1\otimes E)(1\otimes X_t^*)$ moves from $U$ into $U_1\equiv 
U(1\otimes E)(1\otimes Y)U(1\otimes Y^*)+1\otimes (1-E)$. Note 
that the length of this path is at most $\pi$. 

Let $W= U(1\otimes E)(1\otimes Y^*)U(1\otimes Y)$, which is a 
unitary in $C(\T)\otimes B(E\Hil)$. Since $1-E$ is of infinite 
rank, letting $\Hil_k=E\Hil$ for all $k\in\Z$, we identify $\Hil$ 
with $\bigoplus_{k\in {\bf Z}}\Hil_k$ and $E\Hil$ with $\Hil_0$. 
Thus we regard $U_1(s)$ as 
 $$
 \cdots\oplus 1\oplus1\oplus W(s)\oplus 1\oplus \cdots,
 $$                                                 
where $W(s)$ is on $\Hil_0$. Let $W_r=\bigoplus_{k\leq 0}1\oplus 
\bigoplus_{k\geq 1}W\in C(\T)\otimes B(\oplus_{k}\Hil_k)$. Let $H$ 
be a self-adjoint operator on $\bigoplus_k\Hil_k$ such that 
$\Spec(H)=[-\pi,\pi]$, and $e^{iH}$ induces the shift to right. 
Then the map $[0,1]\ni t\mapsto U_1 W_r(1\otimes  
e^{-itH})W_r^*(1\otimes e^{itH})$ moves from $U_1$ into $1$. The 
length of this path is at most $2\pi$.                      

Combining these two paths we get the desired one $(V_t,\ 
t\in[0,1])$. Since $\|W(s_1)-W(s_2)\|\leq 2\|U(s_1)-U(s_2)\|$, we 
get $\|V_t(s_1)-V_t(s_2)\|\leq 6\|U(s_1)-U(s_2)\|$ for any 
$t\in[0,1]$. This concludes the proof. 
\end{pf}  

\begin{lem} \label{U4}
For any $\eps>0$ there exists a $\delta>0$ satisfying the 
following condition: If $Z\in \U(C(\T)\otimes M(\calJ_n)_\alpha)$ 
satisfies that  $\max_{|t|\leq 1}\|(\id\otimes 
\alpha_t)(Z)-Z\|<\delta$, then there exists a rectifiable  path 
$(Z_s,\ s\in[0,1])$ in $\U(C(\T)\otimes M(\calJ_n)_\alpha)$ such 
that $Z_0=1$, $Z_1=Z$, 
 $$
       \max_{s\in [0,1]}\max_{|t|\leq1}\|(\id\otimes \alpha_t)(Z_s)-Z_s\|<\eps, 
 $$
and the length of $(Z_s,\ s\in[0,1])$ is less than $4\pi+\eps$. 
Furthermore if $Z(0)=1$ at $0\in\T$, then the path is chosen to 
satisfy that $Z_s(0)=1$ for $s\in[0,1]$.    
\end{lem}  
\begin{pf}                                                    
We will prove this result for $C(\T)\otimes M(\calJ_n)_\alpha$ 
mostly following the proof of Lemma \ref{U1} for 
$M(\calJ_n)_\alpha$, by using Lemmas \ref{U2} and \ref{U3}, where 
the counterparts for $M(\calJ_n)_\alpha$ are trivial. We will 
indicate below how to use \ref{U2} and \ref{U3} in the proof of 
Lemma \ref{U1}. 

We will apply \ref{U2} when we construct the unitary path $\zeta$ 
in the proof of \ref{U1}.

More explicitly we define $F_k^{\pm}$ as in the proof of \ref{U1}, 
which entails that $F_k^{\pm}\in C(\T)\otimes M(\calJ_n)_\alpha$ 
and 
 \BE
 F_{k}^+ &\leq& 1\otimes E[(2k+1)\delta_0,(2k+3)\delta_0),\\
 F_{2k}^-&\leq & 1\otimes E[(2k-1)\delta_0,(2k+1)\delta_0).
 \EE  
With $\Hil=E[(2k-1)\delta_0,(2k+1)\delta_0)\Hil_0 $, we then 
construct a unitary path $\zeta^k$ in each $C(\T)\otimes B(\Hil)$, 
where $B(\Hil)$ equals 
 $$
 E[(2k-1)\delta_0,(2k+1)\delta_0)M(\calJ_n)_\alpha E[(2k-1)\delta_0,(2k+1)\delta_0).
 $$ 

For $\Hil=E[(4k+1)\delta_0,(4k+3)\delta_0)\Hil_0$, we have to find 
a path $\zeta^{2k+1}$ in $\U(C(\T)\otimes B(\Hil))$ such that 
$\zeta^{2k+1}_0=1$, 
 $$\Ad\,\zeta^{2k+1}_1(F_{2k}^+)=1\otimes 
 E[(4k+1)\delta_0,(4k+2)\delta_0),
 $$ 
and 
 $$\|\zeta^{2k+1}_{t_1}-\zeta^{2k+1}_{t_2}\|\leq \pi|t_1-t_2|.
 $$ 
We can apply \ref{U2} because $1\otimes 
E[(4k+3)\delta_0-\delta',(4k+3)\delta_0)$ is of infinite rank and 
orthogonal to $F_{2k}^+$ and $1\otimes 
E[(4k+1)\delta_0,(4k+2)\delta_0)$ for a small $\delta'>0$. 

For $\Hil=E[(4k-1)\delta_0,(4k+1)\delta_0)\Hil_0$, we find a path 
$\zeta^{2k}$ in $\U(C(\T)\otimes B(\Hil))$ such that 
$\zeta^{2k}_0=1$, $\Ad\zeta^{2k}(F_{2k}^-)=1\otimes 
E[4k\delta_0,(4k+1)\delta_0)$, and 
$\|\zeta_{t_1}^{2k}-\zeta_{t_2}^{2k}\|\leq \pi|t_1-t_2|$. 
                                
Note that $\T\ni s\mapsto Y_k(s)$ is equi-continuous in $k\in\Z$; 
see the definition of $Y_k^+$ in the proof of \ref{U1}.  Since the 
spectrum of $F_k^+$ is contained in $[0,6\mu']\cup [1-6\mu',1]$ 
and $F_k^+$ is the spectral projection of $Y_k^+$ corresponding to 
$[1-6\mu',1]$, it follows that $\T\ni s\mapsto F_k^+(s)$ is 
equi-continuous in $k$. Similarly we get that $\T\ni s\mapsto 
F_k^-(s)$ is equi-continuous in $k$. By \ref{U2} we can choose 
these paths $(\zeta^k)$ such that $\T\ni s\mapsto \bigoplus 
\zeta_t^k(s)$ is continuous for each $t\in [0,1]$, which implies 
that $\zeta_t\equiv \sum_k\zeta_t^k\in C(\T)\otimes M(\calJ_n)$. 
Since $\zeta_t$ commutes with $E[(2k-1)\delta_0,(2k+1)\delta_0)$, 
it follows that $\zeta_t\in C(\T)\otimes M(\calJ_n)_\alpha$ and 
the $\alpha$-spectrum of $\zeta_t$ is contained in 
$[-2\delta_0,2\delta_0]$. Note that $(\zeta_t,\ t\in[0,1])$ has 
length of at most $\pi$.   

When we construct the path $\eta$, we will use Lemma \ref{U3}. 

More explicitly $\zeta_1V$ commutes with $1\otimes E_k$. We apply 
Lemma \ref{U3} and find a path $\eta^k$ in $\U(C(\T)\otimes 
B((E_{2k}+E_{2k+1})\Hil_0))$ which connects 
$\zeta_1V(E_{2k}+E_{2k+1})$ with $E_{2k}+E_{2k+1}$ for each $k$. 
Since $\|\eta^k_t(s_1)-\eta_t^k(s_2)\|\leq 
6\|\zeta_1(s_1)V(s_1)-\zeta_1(s_2)V(s_2)\|$, we have that $\T\ni 
s\mapsto \eta_t(s)\equiv\sum_k\eta_t^k(s)$ is continuous, i.e., 
$\eta_t\in C(\T)\otimes M(\calJ_n)$. Note that the length of 
$(\eta_t,\ t\in [0,1])$ is at most $3\pi$ and that the 
$\alpha$-spectrum of $\eta_t$ is contained in 
$[-4\delta_0,4\delta_0]$. In this way we get the conclusion. 
\end{pf}          

We recall that the flow $\alpha$ on $\calE_n(\subset 
M(\calJ_n)_\alpha)$ induces a flow on the quotient 
$\calO_n(\subset M(\calJ_n)_\alpha/\calJ_n)$, which we will also 
denote by $\alpha$. By using $Q(s_i)\in\calO_n$ we define a unital 
endomorphism $\lambda$ of $\calO_n$ by 
$\lambda(x)=\sum_{i=1}^nQ(s_i)xQ(s_i)^*$. We know that $\lambda$ 
has the Rohlin property (as a version for a single endomorphism or 
automorphism); see \cite{Ror,K95,K96b}. But we have more: 

\begin{lem}\label{U5}                                         
For any $N\in\N$ and $\eps>0$, there are $n^N$ projections $e_{i}$ 
in $\calO_n$ for $i=0,1,\ldots,n^N-1$ such that $e_i\in 
D(\delta_\alpha)$, 
 \BE
 \sum_{i=0}^{2^N-1}e_{i} &=&1,\\ 
 \|\delta_\alpha(e_i)\|&<&\eps,\\
  \max_{i}\|\lambda(e_{i})-e_{i+1}\|&<&\eps,
 \EE
with $F_{n^N}=F_0$.
\end{lem} 
\begin{pf}
Define an action $\gamma$ of $\T^n$ on $\calO_n$ by 
$\gamma_z(s_k)=z_ks_k$ for $z=(z_1,z_2,\ldots,z_n)\in\T^n$. We 
embed $\T$ into $\T^n$ by $z\mapsto (z,z,\ldots,z)$. Then the 
fixed point algebra of $\calO_n$ under $\gamma|\T$ is the closed 
linear span of $Q(s_Is_J^*), \ |I|=|J|$, where $I,J\in 
\{1,2,\ldots,n\}^*$ and is isomorphic to the UHF algebra of type 
$n^\infty$. We will denote it by $\UHF_n$. Then the restriction of 
$\lambda$ onto $\UHF_n$ is the one-sided shift and is known to 
have the Rohlin property. 

Note that $\calO_n^\gamma\subset \UHF_n$, $\calO_n^\gamma\subset 
\calO_n^\alpha$, and by \cite{K05} that $\lambda|\calO_n^\gamma$ 
has the so-called one-cocycle property. (Since $\calO_n^\gamma$ 
has $n$ distinct characters, $\lambda|\calO_n^\gamma$ cannot have 
the Rohlin property but has an approximate Rohlin property.) In 
particular for any $k\in\N$ we have mutually orthogonal 
projections $f_i^k,\ i=0,\ldots,M-1$ in $\calO_n^\gamma$, with 
$M=n^N$, such that 
 $$
 \|\lambda(f^k_i)-f^k_{i+1}\|<1/k
 $$
with $f_{M}^k=f_0$ and $[f^k_0]=[1]$ in $K_0(\calO_n)=\Z/(n-1)\Z$. 
Let $e^k=\sum_{i=0}^{M-1}f^k_i\in \calO_n^\gamma$, which is an 
$\alpha$-invariant projection. Since $\|\lambda(e^k)-e^k\|<M/k$, 
the sequence $(e^k)$ forms a central sequence in $A=\calO_n$. 
Hence $(e^k)$ defines a projection $E$ in $(A'\cap 
A^\omega_\alpha)^\alpha$, where $\omega$ is an ultrafilter on $\N$ 
and $(A'\cap A^\omega_\alpha)^\alpha$ is the $\alpha$-fixed point 
algebra of  the central sequence algebra divided by the ideal 
vanishing at $\omega$. Since $[1]=[e^k]$ in $K_0(A)=\Z/(n-1)\Z$, 
we have that $[1]=[E]$ in $K_0(A'\cap A^\omega)$ \cite{Kp}. (This 
is not entirely trivial; we use the fact that $K_1(A)=0$ for 
$A=\calO_n$.) 

By \cite{K03,Kp} we have that $K_0((A\cap 
A_\alpha^\omega)^\alpha)=K_0(A'\cap A^\omega)$. Since $1,E\in 
(A'\cap A_\alpha^\omega)^\alpha$, we have a $W\in (A'\cap 
A_\alpha^\omega)^\alpha$ such that $WW^*=E$ and $W^*W=1$. Suppose 
that $(w_k)$ represent $W$; then it follows that 
$\lim_\omega\max_{|t|\leq 1}\|\alpha_t(w_k)-w_k\|=0$, 
$\lim_\omega\|w_kw_k^*-e^k\|=0$, and 
$\lim_\omega\|w_k^*w_k-1\|=0$. Thus we can make an isometry $w$ in 
$A=\calO_n$ from $w_k$ for some $k$ with $\eps/3>1/k$ such that 
$\|\lambda(w)-w\|<\eps/3$, $w\in D(\delta_\alpha)$, 
$\|\delta_\alpha(w)\|<\eps/2$, $ww^*=e^k$, and $w^*w=1$. We can 
take the projections $e_i=w^*f_i^kw,\ i=0,1,\ldots,M-1$ because 
then $\sum_ie_i=1$,  
 $$
 \|\lambda(e_i)-e_{i+1}\|<\eps
 $$
and so on. 
\end{pf}

%We recall that $e_n^0=1-\sum_{k=1}^ns_ks_k^*$ is a minimal 
%projection in $\calJ_n$ and that the unital endomorphism $\lambda$ 
%on $M(\calJ_n)\cap \{e_n^0\}'$ is defined by 
%$\lambda(x)=\sum_{k=1}^ns_kxs_k^*+xe_n^0$. Note that 
%$Q\lambda=\lambda Q$ on $M(\calJ_n)\cap \{e_n^0\}'$, where the 
%latter $\lambda$ is the unital endomorphism of 
%$M(\calJ_n)/\calJ_n$ defined by 
%$\lambda(x)=\sum_{k=1}^nQ(s_k)xQ(s_k)^*$.

\begin{lem}
For any $\eps>0$ there exists a $\delta>0$ satisfying the 
following condition: Let $H_\sigma$ be a self-adjoint operator on 
a Hilbert space $\Hil_\sigma$ for $\sigma=1,2$ such that 
$\dim(\Hil_1)=\dim(\Hil_2)<\infty$. If there is a unitary operator 
$W$ from $\Hil_1$ onto $\Hil_2$ such that 
 $$\max_{|t|\leq1}\|e^{itH_2}We^{-itH_1}-W\|<\delta,
 $$
then it follows that $|\lambda^{(1)}_i-\lambda_i^{(2)}|<\eps$ for 
all $i$, where $(\lambda^{(\sigma)}_i)$ is the increasing sequence 
of eigenvalues of $H_\sigma$ (each repeated as often as its
multiplicity) for $\sigma=1,2$. 
\end{lem} 
\begin{pf}
Let $f$ be a non-negative integrable continuous function on $\R$ 
such that $\hat{f}(0)=1$ and $\supp(\hat{f})\subset [-\eps,\eps]$, 
where $\hat{f}(p)=\int f(t)e^{-ipt}dt$. We choose a sufficiently 
small $\delta>0$ so that $W_f=\int f(t)e^{itH_2}We^{-itH_1}dt$ is 
close to $W$ and so is invertible for $W$ satisfying the condition 
in the lemma. Then it follows that for any finite subset $\F$ of 
the eigenvalues of $H_1$ the number of eigenvalues $\lambda$ of 
$H_2$ with $|\lambda-\mu|<\eps$ for some $\mu\in \F$ is greater 
than or equal to the number of elements of $\F$. (Because if $L$ 
is the linear subspace spanned by eigenvectors for $H_1$ whose 
eigenvalues are in $\F$, then the linear subspace $W_fL$ is 
contained in the linear subspace $L'$ spanned by eigenvectors for 
$H_2$ whose eigenvalues are in $\{\lambda\ |\ \exists \mu\in\F\ 
|\lambda-\mu|<\eps\}$. Since $W_f$ is invertible, the dimension of 
$L'$ must be no less than the dimension of $L$.)   This implies, 
by the matching theorem, that there is a bijection $\phi$ from the 
eigenvalues of $H_1$ to those of $H_2$ such that 
$|\phi(\lambda)-\lambda|<\eps$ for any eigenvalue $\lambda$ of 
$H_1$. This implies the above conclusion. 
\end{pf}    

\begin{lem}\label{U7}
For any $\eps>0$ there exists a $\delta>0$ satisfying the 
following condition: If $v\in C(\T)\otimes 
M(\calJ_n)_\alpha/\calJ_n$ is a unitary such that $v(0)=1$ at 
$0\in\T=[0,1]/\{0,1\}$ and $\max_{|t|\leq1} 
\|(\id\otimes\alpha_t)(v)-v\|<\delta$, then there is a unitary 
$Z\in C(\T)\otimes M(\calJ_n)_\alpha$ such that $Q(Z)=v$, 
$Z(0)=1$, and $\max_{|t|\leq1}\|(\id\otimes 
\alpha_t)(Z)-Z\|<\eps$. 
\end{lem} 
\begin{pf}             
Regarding $v$ as a continuous function from $[0,1]$ into 
$\U(M(\calJ_n)_\alpha/\calJ_n)$, we can find a continuous function 
$V$ from $[0,1]$ into $\U(M(\calJ_n)_\alpha)$ such that $V(0)=1$ 
and $Q(V(s))=v(s)$ for $s\in[0,1]$. Since $V(1)\in 1+\calJ_n$ and 
$\U(1+\calJ_n)$ is connected, we may suppose that $V(1)=1$, i.e., 
we have a unitary $V\in C(\T)\otimes M(\calJ_n)_\alpha$ such that 
$V(0)=1$ and $Q(V)=v$. If 
$\max_{|t|\leq1}\|(\id\otimes\alpha_t)(v)-v\|<\delta$, then we 
have that $\max_{|t|\leq1}\|\alpha_t(V(s))-V(s)+\calJ_n\|<\delta$ 
for all $s\in \T$. Thus we may start with such a unitary $V\in 
C(\T)\otimes M(\calJ_n)_\alpha$ instead of $v\in C(\T)\otimes 
M(\calJ_n)_\alpha/\calJ_n$. 

Let $(p_k)$ be an increasing sequence of projections in $\calJ_n$ 
such that $xp_k\ra x$ for any $x\in\calJ_n$ and 
$\alpha_t(p_k)=p_k$ for all $k$. (Note that 
$\alpha_t=\Ad\,e^{iH_0t}$ on $M(\calJ_n)_\alpha$ and $H_0$ is 
diagonal.)  

Since $\T\times [-1,1]\ni (s,t)\mapsto \alpha_t(V(s))-V(s)\in 
M(\calJ_n)_\alpha$ is norm-continuous, the range 
$\{\alpha_t(V(s))-V(s)\ |\ s\in\T,\ t\in[-1,1]\}$ is compact. 
Suppose that $\|\alpha_t(V(s))-V(s)+\calJ_n\|<\delta$ for all 
$(s,t)\in\T\times [-1,1]$ for a sufficiently small $\delta>0$. 
Then there is a $p=p_k$ such that  
 $$\|(\alpha_t(V(s))-V(s))(1-p)\|<\delta
 $$ 
for all $(s,t)\in\T\times [-1,1]$. 

Let $P(s)=V(s)p V(s)^*$ for $s\in\T$. Then $s\mapsto P(s)$ defines 
a projection $P$ in $C(\T)\otimes \calJ_n$ such that 
 $$
 \max_{|t|\leq1}\|(\id\otimes \alpha_t)(P)-P\|<2\delta.
 $$                                                    
Hence, for any $\eps'>0$ by choosing a sufficiently small 
$\delta>0$, we may suppose that there is a unitary $U\in 
C(\T)\otimes (\calJ_n+\C1)$ and a self-adjoint $b\in C(\T)\otimes 
\calJ_n$ such that $\|U-1\|<\eps'$, $\|b\|<\eps'$, and  
$(\id\otimes\delta_\alpha+\ad\,ib)(UPU^*)=0$. We may suppose that 
$U(0)=1$ and $b(0)=0$ since $P(0)=p$ is $\alpha$-invariant. Let 
$P'=UPU^*=UV(1\otimes p)V^*U^*$, which implies that $P'(0)=p$. 
Since $H_0+b(s)$ commutes with $P'(s)$, we can compute 
$\Spec((H_0+b(s))P'(s))$. Let $E_0=\min_{s\in 
\T}\min\Spec((H_0+b(s))P'(s))$ and let $E_1=\max_{s\in \T} 
\max\Spec((H_0+b(s))P'(s))$. 

We find a projection $q\in \calJ_n$ such that $q\leq 1-p$, 
$[H_0,q]=0$, $\Spec(H_0q)\subset [E_0,E_1]$, 
$\rank(q)\gg\rank(p)$, and the eigenvalues of $H_0q$ are almost 
equally distributed in $[E_0,E_1]$. By the last condition we 
should have, in particular, the following: If 
$(\lambda_k)_{k=1}^{\rank(q)}$ is the increasing sequence of 
eigenvalues of $H_0q$ (each of which is repeated as its 
multiplicity indicates), we get $|\lambda_{k}-\lambda_\ell|<\eps$ 
for any $k$ and $\ell$ with $|k-\ell|\leq \rank(p)$.                           

Let $Q=V(1\otimes q)V^*$, which satisfies that $Q(0)=q$ and
 $$
 \max_{|t|\leq1}\|(\id\otimes \alpha_t)(Q)-Q\|<2\delta.
 $$                                                 
Since we can choose the projection $q$ so that $\|(1\otimes 
q)P'\|\approx0$ and $\|QP'\|\approx0$ as precisely as we like, we 
may suppose that the above $U\in C(\T)\otimes (\calJ_n+\C1)$ and 
$b\in\C(\T)\otimes \calJ_n$ are chosen to satisfy that 
$(1\otimes(\delta_\alpha+\ad\,ib))(UQU^*)=0$ in addition, by 
keeping the original relation $P'=UPU^*$. Note, for each $s\in\T$,  
that $U(s)V(s)q$ is a partial isometry from $q\Hil_0$ onto 
$U(s)V(s)q\Hil_0$ and that 
 $$
 \max_{|t|\leq1}\|e^{it(H_0+b(s))}U(s)V(s)qe^{-itH_0}-U(s)V(s)q\|<
  \delta+3\eps',
 $$             
where we have use that $\|b\|<\eps'$, $\|U-1\|<\eps'$, and
 $$
 \|e^{itH_0}V(s)e^{-itH_0}q-V(s)q\|<\delta,\ \ t\in[-1,1].
 $$
By the previous lemma, if $(\lambda_i)$ is the increasing sequence 
of eigenvalues of $H_0q=H_0 q_0$ and $(\mu_i)$ is the increasing 
sequence of eigenvalues of $(H_0+b(s))q_s$ with 
$q_s=U(s)V(s)qV(s)^*U(s)^*$, it follows, by choosing sufficiently 
small $\delta$ and $\eps'$, that $|\lambda_i-\mu_i|<\eps$ for all 
$i$. 

If $S_r$ denotes the set of eigenvalues of $H_0r=rH_0$ with $r$ a 
projection, then $S_{p+q}=S_p\cup S_q$ (which is a disjoint union 
since we count the multiplicity). We have assumed that the rank of 
$q$ is much larger than the rank of $p$; hence the cardinality of 
$S_q$ is much larger than that of $S_p$. We have also assumed that 
$S_q$ overwhelms $S_p$; if we align $S_{p+q}$ in the increasing 
order, the difference of a pair of values which differ by the rank 
of $p$ in this order is at most $\eps$. 

Note that $S_{p+q}(s)=S_p(s)\cup S_q(s)$ for each $s\in\T$, where 
$S_r(s)$ is the set of eigenvalues of $(H_0+b(s))U(s)V(s)r 
V(s)^*U(s)^*$. Since $S_q(s)$ is close to $S_q=S_q(0)$ (as 
$\|b(s)\|<\eps'\ll \eps$) and $S_q$ overwhelms $S_p(s)$, we can 
conclude the following: If $(\lambda_i)$ is the increasing 
sequence of eigenvalues of $H_0(p+q)$ and $(\lambda^{(s)}_i)$ is 
the increasing sequence of eigenvalues of 
$(H_0+b(s))U(s)V(s)(p+q)V(s)^*U(s)^*$, then 
$|\lambda_i-\lambda_i^{(s)}|<2\eps$ for all $i$.   

There are continuous functions $(f_i)$ on $\T$ such that 
 $$
 f_i(s)\in \Spec((1\otimes H_0+b(s))\Ad(U(s)V(s))(p+q))
 $$
for all $s\in\T$, $f_i\leq f_{i+1}$, and 
$(f_i(s))=(\lambda^{(s)}_i)$ as sequences. By perturbing $1\otimes 
H_0+b$ restricted to $UV(1\otimes (p+q))V^*U^*$ slightly, we may 
suppose that there are projections $p_i\in C(\T)\otimes \calJ_n$ 
such that $(1\otimes H_0+b)(UV(1\otimes(p+q))V^*U^*)=\sum_if_ip_i$ 
(see \cite{BBEK}). By a further perturbation up to $\eps$, we may 
suppose that $f_i(s)=\lambda_i,\ s\in\T$ for all $i$. By using 
this, we can construct a partial isometry $W\in C(\T)\otimes 
\calJ_n$ such that $W^*W=1\otimes (p+q)$, 
$WW^*=UV(1\otimes(p+q))V^*U^*$, and $e^{it(1\otimes 
H_0+b)}We^{-itH_0}=W$ for all $t$, where $b$ has now the estimate 
$\|b\|<2\eps$. Since $U(0)V(0)(p+q)V(0)^*U(0)^*=p+q$, we have that 
$W(0)=p+q$. Then it follows that $U^*W$ is a partial isometry from 
$1\otimes (p+q)$ to $P+Q=V(1\otimes (p+q))V^*$ with 
$(U^*W)(0)=p+q$ and that 
 $$
 \max_{|t|\leq1}\|e^{it(1\otimes H_0)}U^*We^{-it(1\otimes H_0)}-U^*W\|
 \leq 2\|U-1\|+2\|b\|.
 $$  
Thus we can set $Z=V(1-1\otimes (p+q))+U^*W$, which is the desired 
unitary in $C(\T)\otimes M(\calJ_n)_\alpha$. 
\end{pf}     
               
\begin{lem}
If $e$ is a projection in $M(\calJ_n)_\alpha/\calJ_n\cap 
\D(\delta_\alpha)$, there is a projection $E\in 
M(\calJ_n)_\alpha\cap \D(\delta_\alpha)$ such that $Q(E)=e$. 
\end{lem}
\begin{pf}
From the proof of Lemma \ref{C1} we can choose $B\in 
M(\calJ_n)_\alpha\cap \D(\delta_\alpha)$ such that $Q(B)=e$ and 
$B^*=B$. Since $Q(B)$ is a projection, the essential spectrum of 
$B$ is $\{0,1\}$ (except for the trivial cases), i.e., $\Spec(B)$ 
has a gap between $0$ and $1$. We can define a projection $E$ from 
$B$ by a $C^\infty$-functional calculus such that $E\in 
\D(\delta_\alpha)$ and $Q(E)=e$. 
\end{pf} 

When $e_0\in  M(\calJ_n)_\alpha/\calJ_n\cap \D(\delta_\alpha)$ is 
a projection such that $\|\delta_\alpha(e_0)\|$ is small, we set 
$h=-i(\delta_\alpha(e_0)e_0-e_0\delta_\alpha(e_0))\in 
M(\calJ_n)_\alpha/\calJ_n$. There is an $H\in M(\calJ_n)_\alpha$ 
such that $H^*=H$, $Q(H)=h$, and 
$\|H\|=\|h\|=\|\delta_\alpha(e_0)\|$. If $E_0\in 
M(\calJ_n)_\alpha\cap \D(\delta_\alpha)$ is a projection such that 
$Q(E_0)=e_0$ as in the above lemma, then we have that 
$[H_0-H,E_0]\in\calJ_n$. Hence, setting 
$K=-i((\delta_\alpha-\ad\,iH)(E_0)E_0-E_0(\delta_\alpha-\ad\,iH)(E_0))\in 
\calJ_n$, we have that $[H_0-H-K,E_0]=0$. Thus we can talk about 
the essential spectrum of $(H_0-H-K)E_0$, which is independent of 
particular choice of $H$ and $E_0$ and may be obtained as follows. 

We define a homomorphism $\phi$ of $C_0(\R)$ into 
$e_0(M(\calJ_n)_\alpha/\calJ_n)e_0$ by 
 $$
 \phi(f)=Q(f(H_0-H-K)E_0),\ \ f\in C_0(\R).
 $$   
Then the kernel of $\phi$ is identified with $\{f\in C_0(\R)\ |\ 
f|S=0\}$ for some closed subset $S$ of $\R$. The essential 
spectrum of $(H_0-H-K)E_0$ is equal to $S$.   

If the above $e_0$ belongs to $\calO_n$, then $h\in\calO_n$ and 
$H\in \calE_n$. The above map $\phi$ can be regarded as mapping 
$C_0(\R)$ into the crossed product $\calO_n\times_\alpha \R$, 
which is known to be simple. (The \cstar\ generated by $\calE_n$ 
and $f(H_0),\ f\in C_0(\R)$ is a quotient of 
$\calE_n\times_\alpha\R$. The quotient by the maximal ideal
$\calJ_n\times_\alpha\R\cong \calJ_n\otimes C_0(\R)$ is isomorphic 
to $\calO_n\times_\alpha\R$.) More precisely, $\phi$ is a 
homomorphism of $C_0(\R)$ into $e_0(\calO_n\times 
\R)e_0=e_0\calO_n e_0\times_{\alpha'}\R$ sending $f$ into 
$\lambda(f)$, where $\alpha'$ is the flow on $e_0\calO_n e_0$ 
generated by $(\delta_\alpha-\ad\,ih)|e_0\calO_ne_0$ and $\lambda$ 
is the canonical unitary group implementing $\alpha'$ in 
$M(e_0\calO_ne_0\times_{\alpha'}\R)$.  Hence if $e_0$ is a 
non-zero projection in $\calO_n$, the essential spectrum of 
$(H_0-H-K)E_0$ is always $\R$.     

The following is a version of Lemma \ref{U4} for 
$M(\calJ_n)_\alpha/\calJ_n$ in place of $M(\calJ_n)_\alpha$.

\begin{lem} \label{calkin}
For any $\eps>0$ there exists a $\delta>0$ satisfying the 
following condition: If $e\in \calO_n\cap D(\delta_\alpha)\subset 
M(\calJ_n)_\alpha/\calJ_n$ is a projection and  $u\in 
(C(\T)\otimes M(\calJ_n)_\alpha/\calJ_n)\cap D(\id\otimes 
\delta_\alpha)$ is a partial unitary such that  $u^*u=1\otimes 
e=uu^*$, $\|(\id\otimes \delta_\alpha)u\|<\delta/2$, and $u(0)=e$ 
at $0\in\T$, then there is a rectifiable path $(w_s,\ s\in[0,1])$ 
of partial unitaries in $C(\T)\otimes M(\calJ_n)_\alpha/\calJ_n$ 
such that $w_s^*w_s=1\otimes e=w_sw_s^*$, $w_0=1\otimes e$, 
$w_1=u$, $w_s(0)=e$, $\max_{|t|\leq1}\|\alpha_t(w_s)-w_s\|<\eps$, 
and the length of $(w_s)$ is less than $4\pi+\eps$. \end{lem}  
\begin{pf}                                                  
We will decide $\delta>0$ later and suppose that $u$ is given as 
above.

We have that $\|\delta_\alpha(e)\|<\delta$. Let 
$h=-i(\delta_\alpha(e)e-e\delta_\alpha(e))\in 
M(\calJ_n)_\alpha/\calJ_n$, which satisfies that 
$\delta_\alpha(e)=\ad\,ih(e)$ and has norm less than $\delta$. We 
find a self-adjoint $H\in \calE_n\subset M(\calJ_n)_\alpha$ such 
that $Q(H)=h$ and $\|H\|=\|h\|$. 

Let $E\in \calE_n\cap D(\delta_\alpha)$ be a projection such that 
$Q(E)=e$. Note that $(\delta_\alpha-\ad\,iH)(E)\in \calJ_n$. Hence 
there is a self-adjoint $K\in \calJ_n$ such that 
 $$
 (\delta_\alpha-\ad\,i(H+K))(E)=0.
 $$
 
Since $u(0)=e$, we can find a continuous path $(U(s),\ s\in[0,1])$ 
in $\U(E M(\calJ_n)_\alpha E)$ such that $U(0)=E$ and 
$Q(U(s))=u(s)$, where $\T$ is regarded as the quotient space 
$[0,1]/\{0,1\}$. Since $Q(U(1))=e$ or $U(1)\in E(\calJ_n+1)E$ and 
the unitary group of $E(\calJ_n+1)E$ is connected, we may suppose 
that $U(1)=E$, i.e., we have a lifting $u$ to the unitary $U$ in 
$C(\T)\otimes EM(\calJ_n)_\alpha E$ such that $U(0)=E$. Let 
$\beta$ be the flow on $M(\calJ_n)_\alpha$ generated by 
$\delta_\alpha-\ad\,iH-\ad\,iK$, which leaves $E$ invariant, i.e., 
$\beta_t$ is implemented by $e^{it(H_0-H-K)}$. Since $                                 
Q(\beta_t(U(s))-U(s))=\alpha_t^{(-h)}(u(s))-u(s)$, we have that 
 $$
 \max_{|t|\leq1}\|\beta_t(U(s))-U(s)+\calJ_n\|<2\|h\|+\delta<3\delta,
 $$                                                                  
where $\alpha^{(-h)}$ is the flow on $M(\calJ_n)_\alpha$ generated 
by $\delta_\alpha-\ad\,ih$. We have remarked that the essential 
spectrum of $(H_0-H-K)E$ is $\R$, where we may suppose that 
$(H_0-H-K)E$ is diagonal by changing $K$ if necessary. We will 
apply Lemmas \ref{U7} and \ref{U4} below for the flow 
$\beta|M(E\calJ_nE)_\beta$ in place of $\alpha|M(\calJ_n)_\alpha$, 
where $E\calJ_nE$ is the compact operators on $E\Hil_0$ as 
$\calJ_n$ is the compact operators on $\Hil_0$ and 
$M(E\calJ_nE)_\beta=EM(\calJ_n)_\beta E=EM(\calJ_n)_\alpha E$. 
This is possible because the properties we needed for 
$\alpha_t=\Ad\,e^{itH_0}$ on $M(\calJ_n)_\alpha$ in these lemmas 
are that $\Spec(H_0)=\R$ and that $H_0$ is diagonal; the same 
properties hold for $\beta_t=\Ad\,e^{it(H_0-H-K)E}$ on 
$M(E\calJ_nE)_\beta$ as asserted above.

We apply Lemma \ref{U7} to the unitary $U$ in $C(\T)\otimes 
EM(\calJ_n)_\alpha E$ (or rather to $Q(U)$) with the Hilbert space 
$E\Hil_0$ and the flow $\beta|EM(\calJ_n)_\alpha E$. Thus we can 
choose a unitary $V\in C(\T)\otimes E M(\calJ_n)_\alpha E$ such 
that $V(0)=E$, $\|U(s)-V(s)+\calJ_n\|<\eps'$ for $s\in\T$ and 
$\max_{|t|\leq1}\|\beta_t(V)-V\|<\eps'$, where $\eps'$ is an 
arbitrarily small constant depending on $\delta$. Then we can 
invoke \ref{U4} to find a rectifiable path $(W_s,\ s\in[0,1])$ in 
$\U(C(\T)\otimes EM(\calJ_n)_\alpha E)$ such that $W_0=E$, 
$W_1=V$, $W_s(0)=E$, $\max_{|t|\leq1}\|\beta_t(W_s)-W_s\|<\eps$, 
and the length of $(W_s)$ is less than $4\pi+\eps$. Now we can set 
$w_s=Q(W_s)$, which satisfies all the required conditions; in 
particular, since $\max_{|t|\leq 1} 
\|\alpha^{(-h)}_t(w_s)-w_s\|<\eps$, it follows that 
$\max_{|t|\leq1}\|\alpha_t(w_s)-w_s\|<\eps+2\delta$. 
\end{pf}

\begin{lem}\label{U8}
For any $\eps>0$ there exists a $\delta>0$ satisfying the 
following condition: If $V\in C(\T)\otimes M(\calJ_n)_\alpha$ is a 
unitary such that  $V(0)=1$ at $0\in\T$ and 
 $$
 \max_{|t|\leq 1}\|(\id\otimes\alpha_t)(V)-V\|<\delta,     
 $$
then there is a unitary $Z\in C(\T)\otimes M(\calJ_n)_\alpha$ such 
that $Z(0)=1$, $\|Q(V)-Q(Z)\lambda(Q(Z))^*\|<\eps$, and 
 $$
 \max_{|t|\leq1}\|(\id\otimes \alpha_t)(Z)-Z\|<\eps.
 $$
\end{lem}       
\begin{pf}          
Let $N\in \N$ and $M=n^N$. we assume that $M^{-1}<\eps$. By 
\ref{U5} there are projections $e_i,\ i=0,\ldots,M-1$ in 
$\calO_n\cap D(\delta_\alpha)\subset M(\calJ_n)_\alpha/\calJ_n$ 
such that $\sum_ie_i=1$, $\|\delta_\alpha(e_i)\|<\eps'$, and 
$\|\lambda(e_i)-e_i\|<\eps'$ for an arbitrarily small $\eps'>0$. 

For each $s\in\T$ there is an isomorphism $\phi_s$ of 
$\calO_n=C^*(Q(s_1),\ldots,Q(s_n))$ onto 
$C^*(Q(V(s)s_1),\ldots,Q(V(s)s_n))\subset M(\calJ_n)_\alpha$. Note 
that $s\mapsto \phi_s$ is continuous in the sense that $s\mapsto 
\phi_s(x)$ is continuous for each $x\in\calO_n$.  Since 
$\phi_s\lambda|\calO_n=\Ad\,Q(V(s))\,\lambda\phi_s|\calO_n$, the 
projections $f_i(s)=\phi_s(e_i)$ satisfy that $\sum_if_i(s)=1$ and 
$\|\Ad\,Q(V(s))\,\lambda(f_i(s))-f_{i+1}(s)\|<\eps'$. Since 
$\alpha_t(V(s))\approx V(s)$ or $\phi_s\alpha_t|\calO_n\approx 
\alpha_t\phi_s|\calO_n$ depending on $\delta$, we can also 
ascertain $\max_{|t|\leq1}\|\alpha_t(f_i(s))-f_i(s)\|<\eps'$, by 
choosing $\delta>0$ sufficiently small. Note that $f_i:s\mapsto 
f_i(s)$ defines a projection in $C(\T)\otimes M(\calJ_n)_\alpha$.  

We assert that there is a partial isometry $z\in C(\T)\otimes
M(\calJ_n)_\alpha/\calJ_n$ such that $z^*z=1\otimes e_0$, 
$zz^*=f_0$, and $\max_{|t|\leq1}\|(\id \otimes \alpha_t)z-z\|$ is 
small. 

Let $h=-i(\delta_\alpha(e_0)e_0-e_0\delta_\alpha(e_0))\in 
(\calO_n)_{sa}$. Then we have that 
$(\delta_\alpha-\ad\,ih)(e_0)=0$. 
                                 
From the proof of \ref{U5} $e_0$ is defined as $w^*p_0w$, where 
$p_0$ is a projection in $\calO_n^\gamma\subset \UHF_n$ and $w$ is 
an isometry in $\calO_n$  with $ww^*\in \calO_n^\gamma$ and 
$ww^*\geq p_0$. We may suppose that $p_0$ and $ww^*$ are in the 
linear span of $Q(s_Is_J^*),\ |I|=|J|=L$ for a large $L\in\N$. 
Since $\phi_s(s_Is_J^*)=Q(V(s))_Ls_Is_J^* Q(V(s))_L^*$ for $I,J\in 
\{1,2,\ldots,n\}^*$ with $|I|=|J|\leq L$, it follows that 
 $$
 f_0(s)=\phi_s(e_0)=\phi_s(w^*)Q(V(s))_Lp_0Q(V(s))_L^*\phi_s(w),
 $$
where $Q(V(s))_L$ is defined inductively by $Q(V(s))_0=1$ and 
 $$Q(V(s))_L=Q(V(s))\lambda(Q(V(s))_{L-1}).
 $$
Thus we can set
 $$
 z(s)=\phi_s(w^*)Q(V(s))_Lwe_0, \ \  s\in\T,
 $$
which has the required property as a partial isometry from 
$1\otimes e_0$ onto $f_0$ if we choose $\delta$ sufficiently 
small. (For example, 
$z(s)^*z(s)=e_0w^*Q(V(s))_L^*\phi_s(ww^*)Q(V(s))_Lwe_0=e_0w^*\cdot 
ww^*\cdot we_0=e_0$.) Since $V(0)=1$ and $\phi_0=\id$, we have 
that $z(0)=e_0$. (Given $\eps>0$ we choose $N\in\N$ and the 
projections $e_i$ for $i=0,1,\ldots, n^N-1$; in particular the 
above $w$, $p_0$, and $L$. Then making $\delta>0$ sufficiently 
small depending on $w$ and $L$, we get that $\phi_s(w^*)$ and 
$Q(V)_L$ are almost $\alpha$-invariant.)

Now we proceed as follows. We choose unitaries $\zeta,\eta$ in 
$C(\T)\otimes M(\calJ_n)_\alpha/\calJ_n$ such that 
$\zeta(s)=\zeta(0)=\eta(0)$, $\zeta\approx1$, 
$\Ad\zeta\lambda(1\otimes e_i)=1\otimes e_{i+1}$, $\eta\approx1$, 
and $\Ad(\eta Q(V))\lambda(f_i)=f_{i+1}$, where $\lambda$ denotes 
the endomorphism $\id\otimes \lambda$ on $C(\T)\otimes 
M(\calJ_n)_\alpha/\calJ_n$. We define a map $\mu$ on $C(\T)\otimes 
M(\calJ_n)_\alpha/\calJ_n$ by $\mu=L_{\eta 
Q(V)}R_{\zeta^*}\lambda$, where $L_x$ (resp. $R_x$) is the left 
(resp. right) multiplication of $x$. Note that $\mu|C(\T)\otimes 
1=\id$, $\mu(z)=\eta Q(V) \lambda(z)\zeta^*$ is a partial isometry 
from $1\otimes e_1$ to $f_1$, etc. We let 
 $$
 y=z^*\mu^{M}(z)\approx z^*Q(V)\lambda(Q(V))\cdots 
 \lambda^{M-1}(Q(V))\lambda^{M}(z)
 $$
which is a unitary in $C(\T)\otimes 
e_0(M(\calJ_n)_\alpha/\calJ_n)e_0$ such that $(\id\otimes 
\alpha_t)(y)\approx y$ for $t\in [-1,1]$ depending on $\delta>0$. 
Since $V(0)=1$ and $z(0)=e_0$, we get that $y(0)=e_0$.   

By assuming $\delta>0$ sufficiently small we can now invoke 
\ref{calkin}, i.e., we find a rectifiable path $(y_s, \ s\in 
[0,1])$ in $C(\T)\otimes M(\calJ_n)_\alpha/\calJ_n$ such that 
$y_sy_s^*=1\otimes e_0=y_sy_s^*$, $y_0=1\otimes e_0$, $y_1=y$, 
$y_s(0)=e_0$, $\|\alpha_t(y_s)-y_s\|\approx0$ for $t\in [-1,1]$, 
and the length of $(y_s)$ is less than $5\pi$. 

We are now in a familiar situation; we define a unitary $u\in 
C(\T)\otimes M(\calJ_n)_\alpha/\calJ_n$ by 
 $$
 u=\sum_{k=0}^{M-1}\mu^k(zy_k),
 $$
where $(y_k)$ is a sequence of unitaries in $C(\T)\otimes 
e_0(M(\calJ_n)_\alpha/\calJ_n)e_0$ chosen from the above path 
$(y_s)$ as follows: $y_0=y$, $\|y_k-y_{k+1}\|<5\pi/M$, 
$y_k(0)=e_0$, and $y_{M-1}=e_0$. Note that $\mu^M(z)=z y$. Then 
$u$ satisfies that $\|\mu(u)-u\|<5\pi/M$, i.e., 
$Q(V)\lambda(u)\approx u$. It also follows that $u(0)=1$ and 
$(\id\otimes \alpha_t)(u)\approx u$ for $t\in [-1,1]$. Then use 
\ref{U7} to lift $u$ to a unitary $Z$ in $C(\T)\otimes 
M(\calJ_n)_\alpha$ such that $Q(Z)\approx u$, $Z(0)=1$, and 
$(\id\otimes \alpha_t)Z\approx Z$ for $t\in [-1,1]$, which 
satisfies the required conditions. 
\end{pf}

\medskip                    
\noindent {\em Proof of Lemma \ref{C4}} 

Let $\delta>0$, which will be decided during the proof. Suppose 
that we are given a unitary $V\in M(\calJ_n)_\alpha$ such that  
$Q(V)\in (\calO_n)'$ and 
$\max_{|t|\leq1}\|\alpha_t(V)-V\|<\delta$. Then by applying Lemma 
\ref{U1} there is a rectifiable path $(V_s,\ s\in[0,1])$ in 
$\U(M(\calJ_n)_\alpha)$ such that  $V_0=1$, $V_1=V$, and 
 $$
 \eps\equiv\max_{s\in [0,1]}\max_{|t|\leq1}\|\alpha_t(V_s)-V_s\|
 $$                                        
is very small depending on $\delta$.   

Since $Q(V\lambda(V)^*)=1$, the map $s\in[0,1]\mapsto 
Q(V_s)\lambda(Q(V_s)^*)$ defines a unitary $w\in C(\T)\otimes 
M(\calJ_n)_\alpha/\calJ_n$ such that $w(0)=1$ and
$\max_{|t|\leq1}\|(\id\otimes\alpha_t)(w)- w\|\leq 2\eps$ for 
$t\in [-1,1]$. By \ref{U7} we obtain a unitary $W\in C(\T)\otimes 
M(\calJ_n)_\alpha$ such that $Q(W)=w$, $W(0)=1$, and 
$\max_{|t|\leq1}\|(\id\otimes\alpha_t)(W)-W\|\equiv \eps'$ is 
small depending on $\eps$. Then by \ref{U8} we obtain a unitary 
$Z\in C(\T)\otimes M(\calJ_n)_\alpha$ such that $Z(0)=1$, 
$w=Q(W)\approx Q(Z^*)\lambda(Q(Z))$, and $\max_{|t|\leq 
1}\|(\id\otimes\alpha_t)(Z)-Z\|$ is small depending on $\eps'$. 
Set $\tilde{V}_s=Z(s)V_s,\ s\in [0,1]$. Then $\tilde{V}_0=1$, 
$\tilde{V}_1=V$, $\lambda(Q(\tilde{V}_s))\approx Q(\tilde{V}_s)$, 
and 
$\max_{s\in[0,1]}\max_{|t|\leq1}\|\alpha_t(\tilde{V}_s)-\tilde{V}_s\|\approx0$. 
Thus we find the desired path $(\tilde{V}_s,\ s\in [0,1])$ 
connecting $1$ with $V$ in $\U(M(\calJ_n)_\alpha)$. This concludes  
the proof.

\section{Cocycle conjugacy}

When $A$ is a unital \cstar, we denote by $\U(A)$ the unitary 
group of $A$. When $C$ is a C$^*$-subalgebra of $A$ such that 
$A\cap C'=C$, we call $C$ a masa of $A$. When $C$ is a masa of $A$ 
and $u\in\U(A)$, $u$ is said to normalize $C$ if $uCu^*=C$. We 
denote by $\Nr(C)$ the set of those $u\in\U(A)$ normalizing $C$. 
Note that $\Nr(C)\supset \U(C)$ and $\Nr(C)$ is a closed subgroup 
of $\U(A)$. When $\Nr(C)$ generates $A$ as a \cstar\ (or 
equivalently the closed linear span of $\Nr(C)$ is $A$), $C$ is 
called a regular masa of $A$.                      

When a regular masa $C$ of $A$ satisfies the conditions that there 
is a norm-one projection of $A$ onto $C$ and that there is a 
character of $C$ which uniquely extends to a state of $A$, we call 
$C$ a {\em Cartan} masa of $A$. 

When $\alpha$ is a flow on a \cstar\ $A$, we denote by 
$\delta_\alpha$ its generator and by $D(\delta_\alpha)$ the domain  
of $\delta_\alpha$, which is a dense $*$-subalgebra of $A$.  

\begin{theo}  \label{AA}
Let $A$ be a unital simple \cstar\ and let $C$ be a Cartan masa of 
$A$. 

Let $\alpha$ be a flow on $A$ such that $D(\delta_\alpha)\supset 
C$. Then the following conditions are equivalent: 
 \begin{enumerate}
 \item $\sup\{\|\alpha_t\delta_\alpha(x)-\delta_\alpha(x)\|\ |\ x\in C,\|x\|\leq 1\}
        \ra 0$ as $t\ra0$.
  \item $\delta_\alpha|C$ is inner, i.e., there is an $h=h^*\in A$ 
        such that $\delta_\alpha(x)=\ad\,ih(x),\ x\in C$.
 \end{enumerate}
\end{theo}

\begin{rem} \label{AA1}
The first condition of the above theorem follows if 
$D(\delta_\alpha^2)\supset C$. This is because then 
$\delta_\alpha^2|C$ is bounded (see, e.g., \cite{Br86}) and 
 $$
 \|\alpha_t\delta_\alpha(x)-\delta_\alpha(x)\|\leq 
 \|\delta_\alpha^2|C\|\,|t|\,\|x\|
 $$
for $x\in C$. \end{rem}   
                              
In the above theorem it is obvious that (2)$\Rightarrow$(1). 
Because, by using $h$ in (2), we have the estimate that 
 $$
 \|\alpha_t\delta_\alpha(x)-\delta_\alpha(x)\|\leq
 2(\|\alpha_t(h)-h\|\,\|x\|+\|h\|\,\|\delta_\alpha|C\|\,|t|),\ \ x\in C.
 $$

Now we assume (1). To derive (2) we first show the following 
lemmas, the first of which is proved in a more general context. 

\begin{lem} \label{AA2}
Let $\alpha$ be a flow on a \cstar\ $A$ and let $C$ be an abelian 
C$^*$-subalgebra of $A$ such that $C\subset D(\delta_\alpha)$. 
Then for each character $\omega$ of $C$ there is a state $\varphi$ 
of $A$ such that $\|\varphi\alpha_t-\varphi\|\ra0$ as $t\ra0$ and 
$\varphi|C=\omega$. \end{lem}    
\begin{pf}               
Let $\omega$ be a character of $C$. Then there is a decreasing net 
$(z_\nu)$ in $\{x\in C\ |\ \|x\|\leq1\}$ such that 
$\omega(z_\nu)=1$ and $(z_\nu)$ converges to the minimal 
projection $p\in C^{**}$ with $\omega(p)=1$. We regard 
$C^{**}\subset A^{**}$ and set $c=\|\delta_\alpha|C\|$. Since 
$\|\alpha_t(z_\nu)-z_\nu\|\leq c|t|$, we get that 
$\|\alpha^{**}_t(p)-p\|\leq c|t|$, i.e., $t\mapsto 
\alpha^{**}_t(p)$ is norm-continuous. Let $B$ be the 
C$^*$-subalgebra of $A^{**}$ generated by $A$ and 
$\alpha_t^{**}(p),\ t\in\R$. Then $B$ is left invariant under 
$\alpha^{**}$ and the restriction $\beta=\alpha^{**}|B$ is a flow, 
i.e, is strongly continuous. There is a $\beta$-cocycle $u$ such 
that $\Ad\,u_t\beta_t(p)=p$. Hence there is a 
$\Ad\,u\beta$-invariant pure state $\hat{\varphi}$ of $B$ such 
that $\hat{\varphi}(p)=1$. We let $\varphi=\hat{\varphi}|A$. Since 
$\alpha=\beta|A=\alpha^{**}|A$, we have that 
$\|\varphi\alpha_t-\varphi\|\leq 
\|\hat{\varphi}\beta_t-\hat{\varphi}\|$. Since 
$\|\hat{\varphi}\beta_t-\hat{\varphi}\|\leq 
\|(\hat{\varphi}-\hat{\varphi}\Ad\,u_t)\beta_t\|\leq 2\|u_t-1\|$, 
we have that $\|\varphi\alpha_t-\varphi\|\leq 2\|u_t-1\|$. Since 
$\varphi|C=\omega$, this concludes the proof. 
\end{pf}    
                       
\begin{lem}  \label{AA3}
Let $\omega$ be a character of $C$ such that $\omega$ uniquely 
extends to a state $\varphi$ of $A$, which is necessarily a pure 
state. Then $\pi_\varphi(C)''$ is a completely atomic masa of 
$\pi_\varphi(A)''=B(\Hil_\omega)$, i.e., there is a family 
$\{p_i\}$ of one-dimensional projections on $\Hil_\varphi$ such 
that $\pi_\varphi(C)''$ is generated by $p_i$'s. If $\Phi$ is a 
norm one projection of $A$ onto $C$, then 
$\pi_\varphi\Phi(x)=\sum_{i}p_i\pi_\varphi(x)p_i$ for $x\in A$. 
\end{lem}    
\begin{pf}
The support projection $p$ of $\omega$ in $C^{**}$ is minimal in 
$A^{**}$. Hence $\pi_\varphi^{**}(p)\equiv p$ is a one-dimensional 
projection belonging to $\pi_\varphi(C)''$. Since $C$ is regular, 
the supremum of $\pi_\varphi(u)p\pi_\varphi(u)^*\in 
\pi_\varphi(C)''$ over $u\in\Nr(C)$ should commute with 
$\pi_\varphi(A)''$; hence equals $1$. Thus we conclude that there 
is a family $\{p_i\}$ of one-dimensional projections in 
$\pi_\varphi(C)''$ such that $\sum_ip_i=1$. 

Let $x\in A$. If $z\in C$ satisfies that $0\leq z\leq 1$ and 
$\pi_\varphi(z)\Omega_\varphi=\Omega_\varphi$, then 
 $$
 \pi_\varphi\Phi(x)\Omega_\varphi=\pi_\varphi\Phi(x)\pi_\varphi(z^2)\Omega_\varphi
 =\pi_\varphi\Phi(zxz)\Omega_\varphi.
 $$
Since $zxz\approx \varphi(x)z^2$ for a suitably chosen $z$, we get 
that 
 $$
 \pi_\varphi\Phi(x)\Omega_\varphi=\varphi(x)\Omega_\varphi=p\pi_\varphi(x)p\Omega_\varphi
 $$
for $x\in A$, where $p$ is the projection onto $\C\Omega_\varphi$. 
Let $u\in \Nr(C)$. Since $x\in C\mapsto 
\lan\pi_\varphi(x)\pi_\varphi(u)\Omega_\varphi,\pi_\varphi(u)\Omega_\varphi\ran= 
\omega(u^*xu)$ is also a character, we apply the same argument and 
get that if $u\in \Nr(C)$, 
 $$
 \pi_\varphi\Phi(x)\pi_\varphi(u)\Omega_\varphi=\varphi(u^*xu)\pi_\varphi(u)\Omega_\varphi
 =p_u\pi_\varphi(x)p_u\pi_\varphi(u)\Omega_\varphi
 =(\sum_{i}p_i\pi_\varphi(x)p_i)\pi_\varphi(u)\Omega_\varphi,
 $$
for $x\in A$, where $p_u=\pi_\varphi(u)p\pi_\varphi(u)^*\in 
\pi_\varphi(C)''$ is the projection onto 
$\C\pi_\varphi(u)\Omega_\varphi$ and $\{p_i\}$ is chosen in the 
first paragraph.  This concludes the proof. 
\end{pf}

\begin{lem} \label{AA4}
Under the condition (1) of the theorem, it follows that 
$\Nr(C)\subset D(\delta_\alpha)$. \end{lem} 
\begin{pf}
Let $\omega$ be a character of $C$ which uniquely extends to a 
state, say $\varphi$, of $A$. Note that $\varphi$ is a pure state 
of $A$.  From \ref{AA2} it follows that 
$\|\varphi\alpha_t-\varphi\|\ra0$ as $t\ra0$. Hence there is a 
unitary flow $U$ on $\Hil_\varphi$ such that 
$\Ad\,U_t\pi_\varphi(x)=\pi_\varphi(\alpha_t(x)),\ x\in A$. We 
will let $H$ be the self-adjoint generator of $U$; $U_t=e^{itH}$. 

Let $h$ be an invariant mean of the function on the abelian group 
$\U(C)$ defined by $u\mapsto -i\pi_\varphi(\delta_\alpha(u)u^*)\in 
B(\Hil_\varphi)$. Then it follows that $h^*=h\in B(\Hil_\varphi)$ 
and $\ad\,ih(\pi_\varphi(x))=\pi_\varphi(\delta_\alpha(x)),\ x\in 
C$. We will show that $t\mapsto \Ad\,U_t(h)$ is norm-continuous. 

Let $c=\|\delta_\alpha|C\|$. We have that 
 $$
 \|\alpha_t(\delta_\alpha(u)u^*)-\delta_\alpha(u)u^*\|
 \leq \|\delta_\alpha(\alpha_t(u)-u)\|+c^2|t|
 $$
for $u\in\U(C)$. Since $\Ad\,U_t(ih)-ih$ is the invariant mean of 
$u\mapsto 
\pi_\varphi(\alpha_t(\delta_\alpha(u)u^*)-\delta_\alpha(u)u^*)$, 
it follows that 
 $$
 \|\Ad\,U_t(h)-h\|\leq \sup_{u\in\U(C)}\|\delta_\alpha(\alpha_t(u)-u)\|+
 c^2|t|.
 $$
Thus we get that $t\mapsto \Ad\,U_t(h)$ is norm-continuous. 

Let $B$ be the \cstar\ generated by $\pi_\varphi(A)$ and 
$\Ad\,U_t(h),\ t\in\R$. Then $\Ad\,U_t$ leaves $B$ invariant and 
$x\mapsto \Ad\,U_t(x)$ is norm-continuous for $x\in B$. Hence 
$t\mapsto \Ad\,U_t|B$ defines a flow on $B$, which we denote by 
$\beta$. Since $h\in B$, $\Ad\,e^{it(H-h)}$ leaves $B$ invariant 
and defines a perturbed flow $\beta^{(-h)}$, which fixes each 
element of $\pi_\varphi(C)$, i.e., $e^{it(H-h)}\in 
\pi_\varphi(C)'$. 

Let $u\in\Nr(C)$ and define 
 $$
 W_t=\pi_\varphi(u)\beta_t^{(-h)}(\pi_\varphi(u^*)),\ \ t\in\R.
 $$                                                             
Then $W_t$ is a unitary in $B$ and $t\mapsto W_t$ is a 
$\beta^{(-h)}$-cocycle. For $z\in\U(C)$ we have that 
$[\pi_\varphi(z),W_t]=0$, which implies that $W_t\in B\cap 
\pi_\varphi(C)'\subset \pi_\varphi(C)''$. Since 
$\beta^{(-h)}_t|B\cap \pi_\varphi(C)'$ is the identity map, we get 
that $t\mapsto W_t$ is a norm-continuous unitary flow, i.e., there 
is a $k=k^*\in B\cap \pi_\varphi(C)''$ such that $W_t=e^{ikt}$. 
Thus it follows that $\pi_\varphi(u)$ is in the domain of the 
generator of $\beta^{(-h)}$ or $\beta$. Since 
$\beta_t\pi_\varphi(u)=\pi_\varphi\alpha_t(u)$ and $\pi_\varphi$ 
is isometric, we get that $u\in D(\delta_\alpha)$. This completes 
the proof of the inclusion $\Nr(C)\subset D(\delta_\alpha)$. 
\end{pf}                                                
        
Continued from the above proof we get that for $u\in\Nr(C)$, 
 $$
 \pi_\varphi(u\delta_\alpha(u^*))-\pi_\varphi(u)\ad\,ih(\pi_\varphi(u^*))
 =ik,
 $$  
where $k=k^*\in \pi_\varphi(C)''$. 

Let $\ol{\Phi}$ be the norm-one projection of $B(\Hil_\varphi)$ 
onto $\pi_\varphi(C)''$. For $z\in C$ we have that 
 $$
 \pi_\varphi\Phi\delta_\alpha(z) =\ol{\Phi}\pi_\varphi\delta_\alpha(z)
 =[i\ol{\Phi}(h),\pi_\varphi(z)]=0.
 $$
Since $h$ is the invariant mean of 
$-i\pi_\varphi(\delta_\alpha(z)z^*),\ z\in \U(C)$ and 
 $$
 \ol{\Phi}(-i\pi_\varphi(\delta_\alpha(z)z^*))
 =-i\pi_\varphi(\Phi(\delta_\alpha(z))z^*)=0,
 $$
we obtain that $\ol{\Phi}(h)=0$. Using that 
$\ol{\Phi}(\pi_\varphi(u)h\pi_\varphi(u^*))=\pi_\varphi(u)\ol{\Phi}(h) 
\pi_\varphi(u^*)=0$ for $u\in \Nr(C)$, we get, from the equation 
in the previous paragraph, 
 $$
 \pi_\varphi(\Phi(u\delta_\alpha(u^*)))=ik,
 $$                                        
i.e., $k\in \pi_\varphi(C)$. Hence we get that 
 $$
 \ad\,ih(\pi_\varphi(u^*))=\pi_\varphi(\delta_\alpha(u^*))-\pi_\varphi(u^*)ik
 \in \pi_\varphi(A).
 $$                
Since this is the case for all $u\in \Nr(C)$, we can conclude that 
$\ad\,ih$ defines a bounded derivation on $\pi_\varphi(A)$. Since 
$\pi_\varphi(A)$ is unital and simple and 
$\pi_\varphi(A)''=B(\Hil_\varphi)$, this implies that $h\in 
\pi_\varphi(A)$ (because any derivation of a unital simple \cstar\ 
is inner). Since 
$\ad\,ih\,\pi_\varphi(z)=\pi_\varphi(\delta_\alpha(z)),\ z\in C$, 
this concludes the proof of Theorem \ref{AA}. 
                                     
\medskip              

In the above theorem, the regularity of the masa $C$ in that 
strong sense is not really needed. We call a masa $C$ of $A$ 
weakly regular if  
 $$ 
 \{u\in \mathcal{PI}(A)\ |\ uu^*,u^*u\in C,\ 
 uCu^*=Cuu^*\}
 $$
generates $A$, where $\mathcal{PI}(A)$ is the set of partial 
isometries of $A$. If a weakly regular masa $C$ satisfies the 
conditions that there is a norm-one projection of $A$ onto $C$ and 
that there is a character of $C$ which uniquely extends to a state 
of $A$, then we call $C$ a weak Cartan masa.  We can get the 
following from the above proof straightforwardly. 

\begin{cor}   \label{oinf}
Let $A$ be a unital simple \cstar\ and let $C$ be a weak Cartan 
masa of $A$. 

Let $\alpha$ be a flow on $A$ such that $D(\delta_\alpha)\supset 
C$. Then the following conditions are equivalent.
 \begin{enumerate}
 \item $\sup\{\|\alpha_t\delta_\alpha(x)-\delta_\alpha(x)\|\ |\ x\in C,\ \|x\|\leq 
 1\}\ra0$ as $t\ra0$. 
 \item $\delta_\alpha$ is inner, i.e., there is an $h=h^*\in A$ 
 such that $\delta_\alpha(x)=\ad\,ih(x),\ x\in C$.
 \end{enumerate}
\end{cor} 

\begin{rem} \label{AA5}
Let $A$ be a unital simple AF \cstar. Let $C$ be a masa of $A$. We 
call $C$ a canonical AF masa if there is an increasing sequence 
$(A_n)$ of finite-dimensional C$^*$-subalgebras of $A$ such that 
the union is dense in $A$ and $C$ is generated by $C\cap A_n\cap 
A_{n-1}',\ n=1,2,\ldots$ with $A_0=0$. If $C$ is a canonical AF 
masa of $A$, then the pair $(A,C)$ satisfies the assumptions of 
Theorem \ref{AA}. Moreover if $\alpha$ is a flow on $A$ such that 
$\D(\delta_\alpha)\supset C$, then $\delta_\alpha|C$ is inner 
\cite{K01}. Note also that a canonical AF masa of $A$ is unique up 
to the transform by automorphisms. 
\end{rem}                   

\begin{prop}\label{BB}
Let $B$ be a stable AF \cstar, $C$ a canonical AF masa, and 
$\gamma$ an automorphism of $B$ such that $\gamma(C)=C$ and 
$\gamma$ fixes no non-trivial ideals of $B$. Suppose that 
$\gamma_*(g)$ is strictly smaller than $g$ for any $g\in 
K_0(B)_+$, i.e., if $\phi$ is a positive homomorphism of $K_0(B)$ 
into $\R$ such that $\phi(g)>0$, then $\phi(g)>\phi\gamma_*(g)$. 
Let $e$ be  a non-zero projection in $C$ such that $\gamma(e)\leq 
e$. Let $A=e(B\times _\gamma \Z)e$ and regard $\calC=Ce$ as a 
C$^*$-subalgebra of $A$. Then it follows that 
 \begin{enumerate}             
 \item $A$ is a unital simple \cstar.
 \item There exists a norm-one projection of $A$ onto $\calC$.
 \item There exists a character of $\calC$ which uniquely extends to a state
 of $A$.    
 \end{enumerate}
If $B$ is simple in addition, $\calC$ is a regular masa of $A$. 
\end{prop}      
\begin{pf}     
\ncm{\calN}{{\Nr}} Since the Connes spectrum of $\gamma$ is full 
(otherwise some power of $\gamma$ would be universally weakly 
inner on some non-zero ideal of $B$), it is well-known that 
$B\times_\gamma\Z$ is simple; thus (1) follows.

By using the dual action of $\T$ on $B\times_\gamma\Z$, we get a 
norm-one projection of $A=e(B\times_\gamma\Z)e$ onto $eBe$. We 
also have a norm-one projection of $eBe$ onto $\calC=Ce$. 
Composing them we get the desired norm-one projection from $A$ 
onto $\calC$. Thus (2) follows. 

To prove (3) we assert that there is a decreasing sequence $(p_k)$ 
of non-zero projections in $\calC$ such that 
$\|p_ku\gamma^n(p_k)\|\ra0$ for all $u\in \calN(\calC)\cap eBe$ 
for all $n\in\N$. Then we take a character $\omega$ of $\calC$ 
such that $\omega(p_k)=1$ for all $k$. If $\varphi$ is an 
extension of $\omega$ to a state of $A$, then we get that 
$\varphi(xU^n)=\lim_k\varphi(p_kx\gamma^n(p_k)U^n)=0$ for all 
$x\in eBe$ and all $n=1,2,\ldots$, where $U$ is the canonical 
unitary multiplier of $B\times_\gamma\Z$ and we have used that 
$\calC$ is a regular masa of $eBe$. This implies that $\varphi$ is 
uniquely determined by $\omega$.     

To prove the above assertion on $(p_k)$, let $(u_m,n_m)_m$ be a 
dense sequence in $(\calN(\calC)\cap eBe)\times\N$. Suppose that 
we have chosen a non-zero projection $p\in \calC$ such that 
$pu_m\gamma^{n_m}(p)=0$ for $m<\ell$. Then we have to find a 
non-zero projection $p'\in\calC$ such that $p'\leq p$ and 
$p'u_\ell \gamma^{n_\ell}(p')=0$. Since 
$q=\Ad\,u_\ell\gamma^{n_\ell}(p)$ is a projection in $\calC$ whose 
equivalence class is strictly smaller than $[p]$, we get that 
$p(1-q)\not=0$; thus we may set $p'=p(1-q)$. We can get the desire 
sequence of projections by repeating this procedure. 

Recall that we have set $\Nr(\calC)=\{u\in \U(A)\ |\ u\calC 
u^*=\calC\}$. To prove the last statement we have to show that 
$\Nr(\calC)$ generates $A$. Note that 
$\lim_{k\rightarrow\infty}\phi([\gamma^k(e)])=0$ for any positive 
homomorphism $\phi$ of $K_0(B)$ into $\R$ (otherwise 
$\lim_k\phi\gamma^k$ would define a non-zero $\gamma$-invariant 
positive homomorphism $\psi$ of $K_0(B)$ into $\R$ such that 
$\psi([e])>0$). We choose $k>1$ such that 
$[\gamma^{k-1}(e)]+[\gamma^k(e)]\leq[e]$ and then choose $v,w\in 
\Nr(\calC)\cap eBe$ such that 
 $$
 \gamma^k(e)\perp \Ad\,v\gamma^{k-1}(e),\ \ \ 
 \gamma^{k-1}(e)\perp \Ad\,w\gamma^k(e),
 $$
and 
 $$
 \gamma^k(e)+\Ad\,v\gamma^{k-1}(e)=\gamma^{k-1}(e)+\Ad\,w\gamma^k(e)\equiv 
 f.
 $$
Then it follows that 
 $$
 S^kS^{*k-1}+vS^{k-1}S^{*k}w^*+e-f\in \Nr(\calC),
 $$
where $S=Ue$. Hence by multiplying $\gamma^{k-1}(e)\in \calC$ from 
the right, we get that $S^kS^{*k-1}=S\gamma^{k-1}(e)$ belongs to 
the \cstar\ $C^*(\Nr(\calC))$ generated by $\Nr(\calC)$. For $u\in 
\Nr(\calC)\cap eBe$ we have that 
$\gamma(u)+e-\gamma(e)\in\Nr(\calC)$ and that 
$\gamma(u)S\gamma^{k-1}(e)u^*=Su\gamma^{k-1}(e)u^*$ belongs to 
$C^*(\Nr(\calC))$. Since $eBe$ is simple, it follows that $S\in 
C^*(\Nr(\calC))$. Since $C^*(\Nr(\calC))\supset eBe$ as follows 
easily, this concludes the proof. 
\end{pf}   

\ncm{\K}{{\mathcal{K}}}      

Given flows $\alpha$ and $\beta$, we say that $\alpha$ is an 
approximate cocycle perturbation of $\beta$ if there is a sequence 
$(u_n)$ of $\beta$-cocycles such that 
$\max_{|t|\leq1}\|\alpha_t(x)-\Ad\,u_n(t)\beta_t(x)\|$ converges 
to zero as $n\ra\infty$ for all $x\in A$. If $\alpha$ is an 
approximate cocycle perturbation of $\beta$, we should note that 
the converse does not follow in general.

We recall that the above automorphism $\gamma$ of $B$ may have the 
Rohlin property (see \cite{EK,BK}). In this case we can show:

\begin{lem} \label{BB1}
Let $A$ be a unital nuclear simple purely infinite \cstar\ and 
$\calC=Ce$ be a Cartan masa of $A$ as in Proposition \ref{BB}. 
Suppose also that the automorphism $\gamma$ of $B$ has the Rohlin 
property. 

Let $\alpha$ be a Rohlin flow on $A$ and let $\beta$ be a flow on 
$A$ such that $\beta_t|\calC=\id$. Then $\beta$ is an approximate 
cocycle perturbation of $\alpha$. 
\end{lem} 
\begin{pf}
Let $(B_n)$ be an increasing sequence of finite-dimensional 
C$^*$-subalgebras of $B$ such that $e,\gamma(e)\in B_1$, the 
central support of $\gamma(e)$ in $eB_1e$ is $e$, 
$\gamma^{\pm1}(B_n)\subset B_{n+1}$, and $C$ is generated by 
$C\cap B_n\cap B_{n-1}',\ n=1,2,\ldots$, where $B_0=0$. Let $u\in 
\Nr(\calC)\cap eB_ne$. We have then $\beta_t(u)=e^{itk}u$ for some 
$k=k^*\in \calC$. Since $\Nr(\calC)\cap eB_ne$ generates $eB_ne$, 
we get that $\beta$ fixes $eBe$ and that $D(\delta_\beta)$ 
contains the union of $eB_ne,\ n=1,2,\ldots$.    

For each $n$ there exists a $z_n=z_n^*\in \calC$ such that 
$\delta_\beta|eB_ne=\ad\,iz_n|eB_ne$. We make a bounded 
perturbation to $\delta_\beta$ and pass to a subsequence of 
$(B_n)$ so that $z_n\in Ce\cap eB_{n+1}e$, i.e., the 
C$^*$-subalgebra $D_n$ of $eB_{n+1}e$ generated by $eB_ne$ and 
$e(B_{n+1}\cap B_n')e\cap C$ is left invariant under $\beta$. Note 
that $eB_ne\subset D_n\subset eB_{n+1}e$.                   

The domain $D(\delta_\alpha)$ may not contain $\bigcup_neB_ne$, 
but we find a $u\in \U(A)$ such that $D(\delta_\alpha)\supset 
ueB_neu^*$ for all $n$. Thus by replacing $\alpha_t$ by a cocycle 
perturbation $\Ad(u^*\alpha_t(u))\alpha_t$, we may suppose that 
$D(\delta_\alpha)\supset \bigcup_neB_ne=\bigcup_nD_n$. 

There exists $y_n=y_n^*\in A$ such that 
$\delta_\alpha|D_n=\ad\,iy_n|D_n$. Then it follows that 
$(\delta_\alpha+\ad i(z_{n}-y_{n}))|D_n=\delta_\beta|D_n$, which 
implies that $\alpha^{(z_{n}-y_{n})}_t(x)=\beta_t(x)$ for $x\in 
D_n$. Thus, for any $n$, by a bounded perturbation on $\alpha$ we 
can always assume that $\alpha=\beta$ on $D_n$. 

Let $S=Ue$ and define $w_t=S^*\beta_t(S)$, which is an 
$\beta$-cocycle. We know that $w_t\in \calC$ and hence that there 
is a $k=k^*\in \calC$ such that $w_t=e^{itk}$. We approximate $k$ 
by $k_n=k_n^*\in Ce\cap eB_ne=Ce\cap D_{n-1}$ such that 
$\|k-k_n\|\ra0$ as $n\ra\infty$.                  

Let $T$ be a large constant and $n$ be a large integer such that 
$1/T\approx0$ and $T\|k-k_n\|\approx0$. Let $N$ be also a large 
integer. We suppose that $\alpha_t(x)=\beta_t(x),\ x\in D_{n+N}$. 

Let $v_t=S^*\alpha_t(S)$, which is a unitary since 
$\alpha_t(SS^*)=SS^*$ and forms a $\alpha$-cocycle. Let $x\in 
Ce\cap eB_{n+N}e=Ce\cap D_{n+N-1}$. Then since $SxS^*\in Ce\cap 
B_{n+N+1}\subset D_{n+N}$, we have that 
$xv_t=S^*SxS^*\alpha_t(S)=S^*\alpha_t(SxS^*S)= v_tx$, which 
implies, in particular, that $e^{-itk_n}v_t=v_te^{-itk_n}$ (as 
$k_n\in Ce\cap D_{n-1}$). Since $\alpha_t(e^{-itk_n})=e^{-itk_n}$, 
it follows that $t\mapsto e^{-itk_n}v_t$ is an $\alpha$-cocycle. 

Let $x\in D_{n+N-1}$. Then since 
$S\beta_{-t}(x)S^*=S\alpha_{-t}(x)S^*\in D_{n+N}$, we have that 
$x\beta_t(S^*)\alpha_t(S)=\beta_t(S^*S\beta_{-t}(x)S^*)\alpha_t(S) 
=\beta_t(S^*)\alpha_t(S)x$, i.e., $e^{-itk}v_t\in A\cap 
D_{n+N-1}'$. Hence, for $x\in D_{n+N-1}$ and $t\in [0,T]$, we get 
that 
$\|[e^{-itk_n}v_t,x]\|=\|[e^{-itk_n}e^{itk},x]e^{-itk}v_t\|\leq 2
T\|k-k_n\|\|x\|$.

Since $\alpha$ has the Rohlin property and $t\mapsto 
e^{-itk_n}v_t$ is an $\alpha$-cocycle, we obtain a unitary $w$ in 
$\U(A)$ such that $\|w\alpha_t(w^*)-e^{-itk_n}v_t\|$ is of the 
order of $|t|/T$. Specifically let us assume that 
$\|w\alpha_t(w^*)-e^{-itk_n}v_t\|\leq \eps$ for $t\in[-1,1]$ with 
$\eps\approx 1/T$. To define such a unitary $w$ we use the path 
$(e^{-itk_n}v_t)_{0\leq t\leq T}$ which connects $1$ with 
$e^{-iTk_n}v_T$, along which any element of the unit ball of 
$D_{n+N-1}$ almost commutes up to the order of $T\|k-k_n\|$, which 
insures that $\|[w,x]\|\leq \eps'\|x\|$ for $x\in D_{n+N-1}$ with 
$\eps'$ of the order of $T\|k-k_n\|$. Since $\|[w,e^{-itk_n}\|\leq 
\eps'$, we get that 
 $$
 \|\alpha_t(Sw)-Swe^{itk_n}\|=\|Sv_t\alpha_t(w)e^{-itk_n}w^*-S\|
 $$                                                                
is less than $\eps+\eps'$ for $t\in[-1,1]$. From the construction 
of $w$, we should also note that $w$ is connected to $1$ by a path 
which almost commutes with $D_{n+N-1}$. We may just as well assume 
that $w\in A\cap D_{n+N-1}'$ with 
$\max_{|t|\leq1}\|\alpha_t(Sw)-Swe^{itk_n}\|\leq\eps+3\eps'$. 

Since $\gamma$ has the Rohlin property and satisfies that 
$\gamma^{\pm}(D_k)\subset D_{k+1}$, we can choose a unitary $u\in 
A\cap (eB_2e)'$ such that $\|w- u\lambda(u^*)\|$ is of order of 
$1/N$ and $u$ belongs to $A\cap D_{n}'$, where $\lambda$ is a 
unital homomorphism of $A\cap (eB_2)'e$ into $A\cap (eB_1e)'$ such 
that $\lambda(x)\gamma(e)=SxS^*,\ x\in A\cap (eB_2e)'$. (To 
construct $u$ we use a path of unitaries  which connects $1$ with 
$w\lambda(w)\lambda^2(u)\cdots \lambda^{M}(w)$ with $M\approx N$ 
and lies in the commutant of $D_n$. This is why we get $u$ from 
$A\cap D_n'$ while $w\in A\cap D_{n+N-1}$; see, e.g., \cite{EK}.) 
Since $\lambda(u)S\lambda(u)^*=Su\lambda(u)^*\approx Sw$, 
$\alpha_t(\lambda(u)S\lambda(u)^*)\approx Swe^{itk_n}$, and 
$[e^{itk_n},\lambda(u)]=0$, we have for $v=\lambda(u)$ that 
 $$
 \Ad\,v^*\alpha_t\Ad\,v(S)\approx v^*Swe^{itk_n}v
 =Su^*we^{itk_n}\lambda(u)\approx Se^{itk_n}\approx \beta_t(S).
 $$
Moreover we have that for $x\in D_{n-1}$ 
 $$
 \Ad\,v^*\alpha_t\Ad\,v(x)=\alpha_t(x)=\beta_t(x).
 $$
This concludes the proof.  
\end{pf}        

When $\alpha$ is a flow on $A$, $\delta_\alpha$ is the generator 
of $\alpha$, and $C$ is a C$^*$-subalgebra of $A$, we say that 
$\alpha$ is $C^{1+\eps}$ if $D(\delta_\alpha)\supset C$ and 
 $$
 \sup_{x\in C,\ \|x\|\leq1}\|(\alpha_t-\id)\delta_\alpha(x)\|
 $$
converges to $0$ as $t\ra0$. 

We recall the following result from \cite{K02}: If $A$ is a unital 
separable nuclear purely infinite simple \cstar, and if each of 
two Rohlin flows on $A$ is an approximate cocycle perturbation of 
the other, then they are cocycle conjugate with each other. By 
\ref{AA} and \ref{BB1} we obtain:

\begin{cor}\label{BB2}
Let $A$ be a unital separable nuclear simple purely infinite 
\cstar\ and $\calC=Ce$ be a Cartan masa of $A$ as in Proposition 
\ref{BB}. Suppose also that $\gamma$ has the Rohlin property. 

Let $\alpha$ and $\beta$ be Rohlin flows on $A$ such that both 
$\alpha$ and $\beta$ are $C^{1+\eps}$ on $\calC$. Then $\alpha$ 
and $\beta$ are cocycle-conjugate. 
\end{cor} 
                          
\begin{pf} By \ref{AA} we may suppose that 
$\alpha|\calC=\id=\beta|\calC$ by inner perturbation. Then by 
\ref{BB1} $\alpha$ is a cocycle perturbation of $\beta$ and vice 
versa. By using the result of \cite{K02} quoted above, we get the 
conclusion. \end{pf} 

Let $\K$ denote the \cstar\ of compact operators on an 
infinite-dimensional separable Hilbert space and let 
$\tilde{\K}=\K+\C1$. We fix a minimal projection $p$ in $\K$. For 
a bounded interval $I$ of $\Z$ we define 
$\tilde{\K}(I)=\bigotimes_I\tilde{\K}$ and 
$\K(I)=\bigotimes_I\K\subset \tilde{\K}(I)$. When $I=[a,b]$, we 
embed $\tilde{\K}([a,b])$ into $\tilde{\K}([a,b+1])$ by $x\mapsto 
x\otimes 1$ and $\tilde{\K}([a,b])$ into $\tilde{\K}([a-1,b])$ by 
$x\mapsto p\otimes x$. Let $\tilde{B}$ denote the inductive limit 
of the system $(\tilde{\K}(I))_I$ with these embeddings and let 
$B$ be the C$^*$-subalgebra of $\tilde{B}$ generated by 
$(\K(I))_I$. Let $\sigma$ denote the  automorphism of $\tilde{B}$ 
induced by the shift on $\Z$ to right; in particular $1\in 
\tilde{\K}(\{0\})$ maps into $p\in \tilde{\K}(\{0\})$. We will 
denote by $e\in \tilde{B}$ the projection corresponding to $1\in 
\tilde{\K}(\{0\})$ or equivalently to $p\in \K(\{-1\})$.      

We note that $\sigma$ leaves $B$ invariant and $e\in B$.
                                                         
Note that $B$ is an AF algebra. It is known \cite{Cun1} that the 
crossed product $B\times_\sigma\Z$ is isomorphic to $\K\otimes 
\calO_\infty$ and that $e(B\times_\sigma\Z)e\cong \calO_\infty$. 

Let $(e_{ij})_{i,j=1}^\infty$ be a family of matrix units in $\K$ 
such that $e_{1,1}=p$ and $\K$ is the closed linear span of 
$(e_{ij})$.  Let $C$ be the abelian C$^*$-subalgebra of 
$\tilde{\K}$ generated by those $e_{i,i}$'s and $1$. It follows 
that $C$ is a regular masa of $\tilde{\K}$. 

Let $C_\infty$ be the C$^*$-subalgebra generated by ${\K}(I)\cap 
\bigotimes_IC$ with all bounded intervals $I$. It follows that 
$e\in C_\infty$ and and that $C_\infty$ is a regular masa of $B$ 
and is left invariant under $\sigma$. Let $U$ denote the unitary 
multiplier of $B\times _\sigma\Z$ which implements $\sigma$. We 
let $S=Ue\in e(B\times_\sigma\Z)e$, which is an isometry. For each 
$x\in\tilde{\K}(\{0\})$ we denote by the same symbol $x$ the 
corresponding element in $B$. 

Under the isomorphism $e(B\times_\sigma\Z)e\cong \calO_\infty$ 
which sends $e_{k1}S$ onto $s_k$, $eC_\infty$ is the abelian 
C$^*$-subalgebra $\calC_\infty$ generated by $s_Is_I^*,\ 
I\in\N^*$, where $\N^*$ is the set of finite sequences in $\N$. It 
is immediate that $\calC_\infty$ is a weak Cartan masa of 
$\calO_\infty$. 

Hence Corollary \ref{oinf} is applicable to the pair 
$A=\calO_\infty$ and $\calC_\infty=eC_\infty$. Since $B$ is not 
simple and $\sigma$ does not have the Rohlin property, we cannot 
apply Lemma \ref{BB1} to this pair $(\calO_\infty,\calC_\infty)$. 
(We have used that $B$ is simple at the beginning of the proof, 
but this is required to define the unital {\em partial} 
endomorphism $\lambda$ of $eBe$ out of $\gamma$, whose Rohlin 
property we would need.) But we can use instead the fact that any 
unital endomorphism of $\calO_\infty$ is approximately inner 
\cite{LP}, in the proof of \ref{BB1}. 

More precisely we define a unital endomorphism $\phi$ of 
$\calO_\infty$ by $\phi(s_k)=s_kw$, or $\phi(e_{k1}S)=e_{k1}Sw$, 
for $k=1,2,\ldots$, where $w$ is the unitary described in the 
proof of \ref{BB1}; in particular $w\in \calO_\infty\cap 
D_{n+N-1}'$. Note that 
$\phi(e_{ij})=\phi(s_is_j^*)=s_iww^*s_j^*=s_is_j^*=e_{ij}$. For 
$e_{i_1,j_1}\otimes e_{i_2,j_2}\in {\K}([0,1])$, if $e_{i_2,j_2}$ 
commutes with $w$, we have that 
 $$
 \phi(e_{i_1,j_1}\otimes e_{i_2,j_2})=\phi(s_{i_1}s_{i_2}s_{j_2}^*s_{j_1}^*)=s_{i_1}w 
 e_{i_2,j_2}w^*s_{j_1}^*=e_{i_1,j_1}\otimes e_{i_2,j_2}.
 $$
In this way, we can show that if $e_k=\sum_{i=1}^ke_{ii}$ is a 
projection in $\K$ such that $w$ is in the commutant of  $D(k)=\C 
e+e_k\K e_k+e_k\K e_k\otimes e_k\K e_k+\cdots+ e_k\K e_k\otimes 
e_k\K e_k\otimes \cdots \otimes e_k\K e_k$ ($k+1$ terms), then 
$\phi|D(k)=\id$. Hence we can find a unitary $v\in \calO_\infty$ 
such that $vSv^*\approx Sw$ and $vxv^*\approx x$ for $x\in D(k)$.  
Since $\bigcup_kD(k)$ is dense in $eBe$, this leads us to the same 
conclusion of \ref{BB1}. 

Thus we have proved Corollary \ref{CC} just as \ref{BB2}. Since 
there are quite a few Rohlin flows on $\calO_\infty$ which are 
trivial on $\calC_\infty$,  this result is certainly non-void. 

We can define a Cartan masa $\calC_n$ of $\calO_n$ in the same way 
as $\calC_\infty$. Since Corollary \ref{BB2} is applicable to the 
pair $(\calO_n,\calC_n)$, let us state: 

\begin{cor} \label{M}               
Let $n=2,3,\ldots$ or $n=\infty$ and let $\calC_n\subset \calO_n$ 
be as above. Then there are Rohlin flows which are $C^{1+\eps}$ on 
$\calC_n$ and  any two Rohlin flows of this type are cocycle 
conjugate with each other. 
\end{cor}                                 

Let $\lambda\in (0,1)$ and let $G_\lambda$ be the subgroup of $\R$ 
generated by $\lambda^n,\ n\in \Z$. Then $G_\lambda$, as an 
ordered subgroup of $\R$, is a simple dimension group and there is 
a stable simple AF algebra $B_\lambda$ such that 
$K_0(B_\lambda)\cong G_\lambda$. Let $\gamma$ be an automorphism 
of $B_\lambda$ such that $\gamma$ induces the multiplication by 
$\lambda$ on $K_0(B_\lambda)=G_\lambda$. We may suppose that 
$\gamma$ leaves a canonical AF masa $C$ of $B_\lambda$ invariant. 
By \cite{EK} $\gamma$ has the Rohlin property and by \cite{Ror} 
the crossed product $B_\lambda\times_\gamma\Z$ is purely infinite 
and simple. (We can get more examples of this kind from 
\cite{BK}.) 
                                                  
If $\{f\in \Z[t]\ |\ f(\lambda)=0\}=p(t)\Z[t]$ for some non-zero 
$p(t)\in \Z[t]$, then $A_\lambda\equiv B_\lambda\times_\gamma\Z$ 
is isomorphic to $\calO_n\otimes \K$ where $n=|p(1)|+1$; otherwise 
$A_\lambda$ is isomorphic to $\calO_\infty\otimes\K$. By cutting 
off $A_\lambda$ by a projection $e\in C$ with $[e]$ a generator of 
$K_0(A_\lambda)$, we get a Cartan masa $\calC=Ce$ in $eA_\lambda 
e$ which is isomorphic to $\calO_n$ with $n$ depending on 
$\lambda$ as above. Thus there are many ways to construct a Cartan 
masa of $\calO_n$ as in \ref{BB}, but we do not know whether we 
can get a new Cartan masa (in case $n<\infty$), which is not 
obtained as an image of the above $\calC_n$ by an automorphism, 
and if we can, whether we have a Rohlin flow which is trivial on 
this Cartan masa.

\medskip
\small
%\begin{flushright}
%Department of Mathematics, Hokkaido University, Sapporo 060 Japan\\
%\end{flushright}

\end{document}